\documentclass[preprint,12pt]{elsarticle}

\usepackage[utf8]{inputenc}
\usepackage{amsmath}
\usepackage{amssymb}
\usepackage[hidelinks]{hyperref}
\hypersetup{
     colorlinks   = true,
     citecolor    = blue,
     linkcolor    = blue}
\usepackage{listings}
\usepackage{xcolor}
\usepackage{geometry}
\usepackage{graphicx}
\usepackage{subcaption}
\usepackage{epic}
\usepackage{comment}
\usepackage{float}

\DeclareSymbolFont{lettersA}{U}{txmia}{m}{it}

\newcommand{\fvec}{\mathbf{f}}
\newcommand{\zvec}{\mathbf{z}}

\definecolor{codegreen}{rgb}{0,0.6,0}
\definecolor{codegray}{rgb}{0.5,0.5,0.5}
\definecolor{codepurple}{rgb}{0.58,0,0.82}
\definecolor{backcolour}{rgb}{0.95,0.95,0.92}

\lstdefinestyle{mystyle}{
    backgroundcolor=\color{backcolour},
    commentstyle=\color{codegreen},
    keywordstyle=\color{magenta},
    numberstyle=\tiny\color{codegray},
    stringstyle=\color{codepurple},
    basicstyle=\ttfamily\footnotesize,
    breakatwhitespace=false,
    breaklines=true,
    captionpos=b,
    keepspaces=true,
    showspaces=false,
    showstringspaces=false,
    showtabs=false,
    tabsize=2
}

\title{Super Time Stepping Methods for Diffusion using Discontinuous-Galerkin Spatial Discretizations}

\author[1,2]{Mustafa Aggul\corref{cor1}}
\author[3]{Manaure Francisquez}
\author[1]{Daniel R. Reynolds}
\author[1]{Sylvia Amihere}

\affiliation[1]{
  organization={Department of Mathematics \& Statistics},
  institution={University of Maryland Baltimore County,},
  city={Baltimore},
  state={MD},
  country={USA}
}

\affiliation[2]{
  organization={Department of Mathematics},
  institution={Southern Methodist University,},
  city={Dallas},
  state={TX},
  country={USA}
}

\affiliation[3]{
  organization={Princeton Plasma Physics Laboratory (PPPL)},
  city={Princeton},
  state={NJ},
  country={USA}
}

\cortext[cor1]{Corresponding author, maggul@umbc.edu}

\begin{document}

\begin{abstract}
Super-time-stepping (STS) methods provide an attractive approach for enabling explicit time integration of parabolic operators, particularly in large-scale, higher-dimensional kinetic simulations where fully implicit schemes are impractical. In this work, we present an explicit STS framework tailored for diffusion operators in gyrokinetic models, motivated by the fact that constructing and storing a Jacobian is often infeasible due to strong nonlocal couplings, high dimensionality, and memory constraints. We investigate the performance of several STS methods, including Runge–Kutta–Chebyshev (RKC) and Runge–Kutta–Legendre (RKL) schemes, applied to a diffusion equation discretized using both discontinuous Galerkin (DG) and finite-difference methods. To support time adaptivity, we introduce a novel error norm designed to more accurately track temporal error arising from DG spatial discretizations, in which degrees of freedom contribute unevenly to the solution error. Finally, we assess the performance of an automatic eigenvalue estimation algorithm for determining the required number of STS stages and compare it against an analytical estimation formula.
\end{abstract}

\begin{keyword}
adaptive super-time-stepping, diffusion, discontinuous Galerkin, dominant eigenvalue estimation
\end{keyword}

\maketitle

\section{Introduction}
\label{sec:intro}

A large number of models based on partial differential equations (PDEs) use parabolic or diffusive terms to represent physical phenomena of interest. Heat transport via diffusion is one such process. In the field of plasma simulation with PDEs, elastic collisions between charged particles are often studied using a model with a diffusive component, as is the case of the Landau-Rosenbluth Fokker-Planck operator~\cite{Rosenbluth1957}. When studying magnetized plasmas for magnetic confinement fusion (MCF) energy, collisions play a central role. Fusion plasmas must retain their energy density in order to produce enough fusion energy for the device to be economically viable, yet turbulent eddies threaten these conditions by transporting energy out of the plasma core towards the vessel walls. This turbulent transport can be mitigated by ``zonal flows", a large-scale self-organized rotation of the plasma (akin to Jupiter's latitudinal flows) which pose barriers for turbulent eddies to escape the plasma core. These zonal flows can be damped by collisions~\cite{Hinton1999}, and thus accurately modeling collisions becomes a necessary ingredient to accurately predict the efficiency of a MCF plasma.

The gold standard model for turbulent transport in magnetized fusion plasmas is the gyrokinetic model. This framework stems from a reduction of the six-dimensional (in position and velocity space) Vlasov kinetic model~\cite{Kardar2007} by averaging over the fast gyromotion of a particle around the magnetic field (in the case of strong magnetization). Such a procedure removes one of the velocity coordinates (the gyro-angle), and orders out time scales that are shorter than the gyro-period, which is appropriate as most of the turbulent transport in MCF is due to modes with frequencies lower than the inverse of the gyro-period. This model prescribes the time evolution of the distribution function $f_s$ of species $s$ (with mass $m_s$ and charge $q_s$) over position ($\mathbf{x}$) and velocity ($v_{\|},\mu$) space using the velocity $v_{\|}=\mathbf{\hat{b}}\cdot\mathbf{v}$ parallel to the magnetic field $\mathbf{B}=B\,\mathbf{\hat{b}}$ and the magnetic moment $\mu=m_sv_{\perp}^2/(2B)$ (here $v_\perp$ is the velocity perpendicular to $\mathbf{B}$) as effective velocity coordinates. Such time evolution is dictated by a PDE such as~\cite{Francisquez_2025} 
\begin{equation}
  \label{eq:gkeyll_model}
  \frac{\partial \mathcal{J} f_s}{\partial t} + \nabla\cdot\left(\mathcal{J} f_s\mathbf{\dot{x}}\right) + \frac{\partial}{\partial v_{\|}}\left(\mathcal{J} f_s\dot{v}_{\|}\right) = 
  \mathcal{J}\mathcal{C}_s.
\end{equation}
The second and third term on the left-hand side provide advection in $\mathbf{x}$ and $v_{\|}$ space, respectively. The term on the right-hand side corresponds to the effect of collisions. A collision model that retains the advective-diffusive structure of the Fokker-Planck operator is the Dougherty operator~\cite{Dougherty1964}, which is used in a variety of gyrokinetic codes~\cite{francisquez_conservative_2020,Jo_2022,Ulbl_2023,Hoffmann_2023,Mandell_2024,Dorf2025}, and consists of
\begin{equation}
  \label{eq:collision_model}
  \mathcal{C}_s = 
  \nu_{ss}\left\{\frac{\partial}{\partial v_{\|}} \left[(v_{\|}-u_{\|})\mathcal{J} f_s + v_{ts}^2 \frac{\partial \mathcal{J} f_s}{\partial v_{\|}}\right] + \frac{\partial}{\partial \mu}\left[2\mu\mathcal{J} f_s + 2\frac{m_sv_{ts}^2}{B}\mu\frac{\partial \mathcal{J} f_s}{\partial \mu}\right]\right\},
\end{equation}
where $u_{\|}=u_{\|}(\mathbf{x},t;f_s)$ and $v_{ts}=v_{ts}(\mathbf{x},t;f_s)$ are the mean flow parallel and thermal speeds, respectively, and $\nu_{ss}$ is the collision frequency. Note that the operator in equation~\eqref{eq:collision_model} is an integro-differential nonlinear operator, since $u_{\|}$ and $v_{ts}$ involve velocity integrals of $f_s$ itself. Near the walls of a fusion device, where plasma density can be moderate and temperature can be very low, the collision frequency can become large ($\nu_{ss}$ is proportional to density and inversely proportional to temperature), and the diffusion terms (proportional to $v_{ts}$), if time-integrated with explicit methods, can pose an impractically challenging stability (CFL) condition that requires a very small time step ($\Delta t$).

Such a severely small $\Delta t$ can be avoided by using implicit methods. However, constructing the matrix arising from treating the diffusive terms in equation~\eqref{eq:collision_model} with a backward Euler method, for example, or for constructing the Jacobian of the operator can be very challenging, in part because as mentioned in the previous paragraph, this operator is nonlocal and nonlinear, resulting in a dense matrix. For this reason, in this work we focus on ``Jacobian-free" methods, where at first we omit the advective terms from equation \eqref{eq:gkeyll_model} and consider only the diffusive term in equation~\eqref{eq:collision_model} -- in future studies we will extend this work into a solver for the full gyrokinetic collision operator in equation~\eqref{eq:collision_model}. We explore temporally-adaptive methods that achieve a requested accuracy tolerance (as opposed to using fixed $\Delta t$ set by heuristic measures such as user experience or linear stability). These methods include diagonally-implicit Runge--Kutta (DIRK) \cite{reynoldsARKODEFlexibleIVP2023}, strong stability preserving (SSP) explicit Runge--Kutta (ERK) \cite{feketeEmbeddedPairsOptimal2022, ketchesonHighlyEfficientStrong2008}, and explicit super time stepping (STS) \cite{sommeijer_rkc_1998, meyerSecondorderAccurateSuper2012, verwerRKCTimesteppingAdvection2004}. A key input for STS methods is a reliable estimate of the dominant eigenvalue of the Jacobian of the diffusion operator in question --  therefore, we explore the performance of an algorithm that automatically estimates this dominant eigenvalue, comparing this against an analytical formula. Furthermore, this work focuses on the application of such time-integration methods to a discontinuous Galerkin (DG) discretization of the problem~\cite{Atkins1998,Francisquez_2025} using the Gkeyll code~\cite{gkeyllWeb}. To that end, we introduce a novel error norm that is designed to balance the relative contributions from solution degrees of freedom to overall solution error, in contrast with standard approaches that weigh each degree of freedom evenly. 

This manuscript is organized as follows. Section 2 presents the time-integration methods used in this work. Section 3 discusses the estimation of the dominant eigenvalue to be used in STS methods. Section 4 provides the results of applying these methods to a representative diffusion problem, and section 5 offers concluding remarks. We note that all test problems in this paper are available in the public GitHub repository \cite{CEDADemonstrationsRepo}, and use time integrators from the v7.5.0 release of SUNDIALS \cite{gardnerEnablingNewFlexibility2022a,roberts2025new}.

\section{Method Descriptions}
\label{sec:method-descriptions}

As stated above, we wish to explore efficient Jacobian-free time-integration of diffusive models posed on DG spatial semi-discretizations.  To that end, we investigate the diffusion problem
\begin{equation}
    \label{eq:diffusion}
    \begin{split}
    \frac{\partial \mathcal{J} f_s}{\partial t}(t,v_{\|},x) &= \frac{\partial}{\partial {v_{\|}}} \left[ D_{v_{\|}}(v_{\|}) \frac{\partial \mathcal{J} f_s}{\partial {v_{\|}}}(t,v_{\|},x)\right],\\
    \mathcal{J} f_s(0,v_{\|},x) &= \mathcal{J} f_{s,0}(v_{\|},x),
\end{split}
\end{equation}
where diffusion is only applied in the $v_{\|}$ direction using a spatially-dependent diffusion coefficient $D_{v_{\|}}(v_{\|})$. As is typical when considering time integration methods, we rewrite the diffusion problem as
a first-order initial-value problem (IVP),
\begin{equation}
  \label{eq:IVP_unsplit}
  \frac{\partial \fvec}{\partial t}(t) = \mathcal{G} (t,\fvec), \quad \fvec(t_0) = \fvec_0,
\end{equation}
where $\mathcal{G}$ encodes the discretized diffusion operator and $\fvec$ is a large vector containing the solution values throughout the finite-difference grid or the DG expansion coefficients throughout the domain. Due to our initial focus on the pure diffusion regime, we assume that the Jacobian $\mathcal{G}_{\fvec}$ has negative, real-valued eigenvalues.  

\subsection{DIRK Methods}

Diagonally-implicit Runge--Kutta methods evolve the equation \eqref{eq:IVP_unsplit} by proceeding through a sequence of stage solutions, $\mathbf{z}_i$, $i=1,\ldots,s$.  At each stage DIRK methods must solve an implicit system of algebraic equations of the form
\begin{equation}
    \label{eq:DIRK_nonlinear}
    \mathbf{0} = \mathcal{F}(\zvec_i) := \zvec_i - h A_{i,i} \mathcal{G}(t_n + c_ih, \zvec_i) - \mathbf{a}_i,
\end{equation}
where $h$ is the time step size, $A_{i,i}$ is the $i$th diagonal entry of the DIRK Butcher tableau and $\mathbf{a}_i$ is a vector of known data (e.g., previous time step, previous stages).  The stage solutions are then combined to construct an approximate solution $\fvec_{n+1}$ and an ``embedded'' solution, $\tilde{\fvec}_{n+1}$, 
\begin{align}
    \label{eq:DIRK_solution}
    \fvec_{n+1} &= \fvec_n + h\sum_{i=1}^s b_i \mathcal{G}(t_n + c_ih, \zvec_i),\\
    \label{eq:DIRK_embedding}
    \tilde{\fvec}_{n+1} &= \fvec_n + h\sum_{i=1}^s \tilde{b}_i \mathcal{G}(t_n + c_ih, \zvec_i),
\end{align}
where the embedding typically has order of accuracy one lower than the solution.  Thus their difference provides an estimate of the local temporal error introduced in the time step,
\begin{equation}
    \label{eq:DIRK_local_error}
    \mathbf{\varepsilon}_{n+1} = \fvec_{n+1} - \tilde{\fvec}_{n+1},
\end{equation}
as will be discussed in Section \ref{sec:TemporalAdaptivity} for performing time step adaptivity.

To solve the equation \eqref{eq:DIRK_nonlinear}, we use a standard inexact Newton method, that uses a preconditioned conjugate gradient (CG) method from SUNDIALS \cite{gardnerEnablingNewFlexibility2022a}.  For the preconditioner
to keep in-line with the goal of solvers that do not require construction of the full Jacobian matrix, we use a simple Jacobi preconditioner, where the preconditioner matrix $P$ is formed using only the diagonal of the Jacobian $\mathcal{F}_{\zvec}$.
In the tests that follow, we use the default embedded DIRK methods of orders 2 and 3 from the SUNDIALS ARKStep solver, corresponding with the implicit portion of the ARK2 method from \cite{giraldo_implicit-explicit_2013}, and the ESDIRK3(2)5L[2]SA method from \cite{kennedy_diagonally_2016}.

We note that although inexact Newton methods do not require the Jacobian matrix $\mathcal{F}_{\zvec}$ to be constructed, our preconditioner requires information about the Jacobian diagonal.  As a result, we only test DIRK methods when performing tests based on a finite-difference spatial semi-discretization, so that we can provide baseline comparisons between adaptive time integration methods in this simplified scenario: The diffusion operator is spatially variable but independent of other solution components, so the implicit solves do not introduce strong global coupling. In contrast, applying DIRK to the full gyrokinetic collision operator (equation \eqref{eq:collision_model}) is generally impractical since each stage requires the solution of a tightly coupled nonlocal system, and hence frequent global communication. Consequently, the per-step computational cost becomes prohibitive for routine use in a production gyrokinetic solver. For tests that use Gkeyll's DG spatial discretization, we therefore omit DIRK methods.
\subsection{SSP ERK Methods}

As with most DG-based codes, the default time integration methods in Gkeyll are Strong-Stability-Preserving (SSP) Runge--Kutta methods \cite{shu1988ssp,gottlieb2009ssp}.
SSP methods are designed to preserve monotonicity and nonlinear stability properties
(e.g., total-variation-diminishing (TVD)) satisfied by the forward Euler method under
a restricted step size. These properties are particularly important for
convection-dominated and hyperbolic flow problems.

Following the Shu--Osher formulation \cite{shu1988ssp}, SSP methods applied to the IVP \eqref{eq:IVP_unsplit} may be expressed as
\begin{align}
\zvec_1 &= \fvec_n, \\
\zvec_i &= \sum_{j=1}^{i-1}
\left( \alpha_{i,j} \zvec_j + \beta_{i,j} h \, \mathcal{G} (t_n + c_j h, \zvec_j) \right),
\qquad i = 2,\ldots,s, \label{eq:ssprk} \\
\fvec_{n+1} &= \zvec_s.
\end{align}
Here $s$ denotes the number of stages, $h$ is the time step, and the coefficients satisfy
\[
\alpha_{i,j} \ge 0, \qquad \beta_{i,j} \ge 0, \qquad
\sum_{j=1}^{i-1} \alpha_{i,j} = 1,
\]
ensuring that each stage is a convex combination of previous stage values and forward
Euler updates.  When the spatial discretization admits monotone stability under forward
Euler with step size $h_{\mathrm{FE}}$, the SSPRK method preserves this property under
\[
h \le C \, h_{\mathrm{FE}},
\]
where $C$ is the SSP coefficient \cite{gottlieb2001strong}.

This Shu--Osher structure enables a memory-efficient low-storage Runge--Kutta (LSRK) implementation, requiring
only two work vectors regardless of the number of stages, $s$. SSPRK methods such as SSP$(s,2)$, SSP$(s,3)$, and
SSP$(10,4)$ \cite{ketchesonHighlyEfficientStrong2008,ketcheson2011runge} are widely used in
large-scale PDE solvers and are included in ARKODE's LSRKStep module within the SUNDIALS suite
\cite{roberts2025new}. In the computations (Section \ref{sec:diffusion-tests}),
the fixed number of stages employed for the second--, third-- and fourth--order SSP methods
are 2, 4 and 10, respectively. These SSP schemes are particularly suitable for high-resolution
finite-volume, finite-difference, and discontinuous Galerkin discretizations for transport-dominated problems, where preserving monotonicity and stability is essential.

We note that like DIRK methods, SSPRK methods can admit embedded solutions $\tilde{\fvec}_{n+1}$ and local temporal error estimates $\mathbf{\varepsilon}_{n+1}$, using the same formulas \eqref{eq:DIRK_embedding} and \eqref{eq:DIRK_local_error}.  
In this work, we use the embedded SSP Runge--Kutta pairs designed by Fekete et al.~\cite{feketeEmbeddedPairsOptimal2022}, that guarantee SSP properties for both the primary
and embedded solutions, ensuring that the embeddings preserve the same monotonicity and stability characteristics as the computed solutions.

\subsection{STS Methods}
\label{sec:STSmethods}

Due to the diffusive nature of the model \eqref{eq:diffusion}, we anticipate that SSP methods will be stability-limited, and that the maximum linearly stable step size will result in solutions with significantly more accuracy than required. Thus, we additionally consider
Super Time Stepping (STS) Runge--Kutta integrators. These explicit methods incorporate additional internal stages to extend the stability region
along the negative real axis, making them suitable for systems with
parabolic stiffness (i.e., Jacobians with large negative real eigenvalues). Unlike fully
implicit approaches, STS methods require no implicit algebraic solves, thereby retaining the simplicity of explicit time-stepping while
providing enhanced stability for stiff diffusion-dominated regimes.

We employ two families of second-order, explicit STS methods that support temporal adaptivity:
the Runge--Kutta--Chebyshev (RKC) scheme \cite{sommeijer_rkc_1998} and Runge--Kutta--Legendre (RKL)
scheme \cite{meyerSecondorderAccurateSuper2012, meyerStabilizedRungeKutta2014a}. These methods
take the form
\begin{subequations}
  \label{eq:RKC2}
  \begin{align}
    \zvec_0 &= \fvec_n,\\
    \zvec_1 &= \zvec_0 + h \tilde{\alpha}_1 \mathcal{G}(t_n, \zvec_0),\\
    \zvec_j &= \alpha_j \zvec_{j-1} + \beta_j \zvec_{j-2} + (1-\alpha_j-\beta_j)\zvec_0 \\
    \notag
    &\quad + h\tilde{\alpha}_j \mathcal{G}(t_{n,j-1}, \zvec_{j-1}) + h\tilde{\gamma}_j \mathcal{G}(t_{n}, \zvec_0), \quad j=2,\ldots,s,\\
    \fvec_{n+1} &= \zvec_s,
  \end{align}
\end{subequations}
The coefficients $\alpha_j$, $\beta_j$, $\widetilde{\alpha}_j$, and $\widetilde{\gamma}_j$ are chosen
to maximize the stability interval while maintaining second-order accuracy; specific values
are provided in \cite{sommeijer_rkc_1998, meyerSecondorderAccurateSuper2012}. The number of internal
stages $s$ typically scales as $\mathcal{O}(\sqrt{\lambda})$, where $\lambda$ is the magnitude of the dominant eigenvalue of the Jacobian $\mathcal{G}_{\fvec}$. This scaling allows STS methods to
dramatically extend the stability limits, enabling significantly larger time steps for
diffusion-dominated problems compared to classical explicit Runge--Kutta methods.

LSRK-type STS methods are particularly effective for large-scale diffusion problems arising
in computational physics, including heat conduction, viscous flows, and plasma transport.
They are implemented in ARKODE's LSRKStep module within the SUNDIALS suite
\cite{roberts2025new}, providing an explicit, memory-efficient alternative to fully implicit
solvers when diffusion stiffness dominates.

Time adaptivity for the RKC and RKL schemes is performed using a local error estimate based on
a cubic Hermite interpolating polynomial. Following the approach in
\cite{sommeijer_rkc_1998, roberts2025new}, the error estimator is defined as
\begin{equation}
\label{eq:RKC-error-estimate}
\mathbf{\varepsilon}_{n+1}
  = \frac{1}{15} \left[
        12\left(\fvec_{n} - \fvec_{n+1}\right)
        + 6h\left(\mathcal{G}(t_n,\fvec_{n})
        + \mathcal{G}(t_{n+1},\fvec_{n+1})\right)
    \right],
\end{equation}
where $\fvec_{n}$ and $\fvec_{n+1}$ denote the numerical solutions at times $t_n$ and $t_{n+1}$,
respectively. This estimator provides a second-order accurate approximation of the local temporal error, enabling adaptive step-size selection.

Given a desired step size $h$, STS methods must estimate the required number of internal stages to ensure linear stability; this in turn requires an estimate of the dominant eigenvalue, $\lambda$.
These estimates can be obtained either through an application-specific function ($\lambda_{user}$), or they may be estimated internally using the new SUNDomEigEst module in SUNDIALS ($\lambda_{approx}$). These estimates are then multiplied by a safety factor to account for any estimation error to ensure that the computed number of stages yields a stable time step. The estimation and the safety factor are discussed in Section \ref{sec:dominant-eigenvalue}.

\section{Temporal Adaptivity}
\label{sec:TemporalAdaptivity}

To control the time steps used by each of the above methods, we leverage the local temporal error estimates $\mathbf{\varepsilon}_{n+1}$ provided by each method.  The first step in this process is to convert the vector-valued $\mathbf{\varepsilon}_{n+1}$ into a scalar-valued measure of ``overall'' temporal error in the step.  By default, SUNDIALS uses a component-wise weighted RMS norm for this purpose:
\begin{equation}\label{eq:norm1}
\|\mathbf{\varepsilon}_{n+1}\|_{WRMS} = \sqrt{\frac{1}{N} \sum_{i=1}^{N}\left(\frac{\mathbf{\varepsilon}_{n+1,i}}{RTOL |\fvec_{n,i}|+ATOL}\right)^2},
\end{equation}
where $N$ denotes the total number of global degrees of freedom, $\fvec_{n,i}$ is the $i$-th component of the vector $\fvec_n$, $RTOL$ corresponds to the desired relative solution error in a step, and $ATOL$ corresponds with a floor below which the absolute error the solution is considered noise \cite{gardnerEnablingNewFlexibility2022a}.
However, a DG discretization may weigh the importance of each degree of freedom differently (e.g. cell average may be weighted more heavily than cell slopes). Thus, the norm \eqref{eq:norm1} may over- or under-estimate individual contributions to the overall solution error. To that end, we also consider a cell-wise norm:
\begin{equation}\label{eq:norm2}
\|\mathbf{\varepsilon}_{n+1}\|_{WRMS} = \sqrt{\frac{1}{N_c}\sum_{i=1}^{N_c}\frac{\|\mathbf{\varepsilon}_{n+1,i}\|_C^2}{\left(RTOL \|\fvec_{n,i}\|_C+ATOL\right)^2}}, \qquad \|\fvec_{n,i}\|_C := \sqrt{\frac{1}{N_b}\sum_{j=1}^{N_b} \fvec_{n,i,j}^2},
\end{equation}
where $N_c$ and $N_b$ refer to the number of cells and the number of degrees of freedom in each cell, respectively, and $\fvec_{n,i,j}$ corresponds with the $j$-th degree of freedom in the $i$-th cell of the vector $\fvec_n$. 

For both choices of norm, \eqref{eq:norm1} and \eqref{eq:norm2}, our temporal adaptivity follows the standard step-size control procedures in ARKODE \cite{reynoldsARKODEFlexibleIVP2023}, that strive to ensure the local error satisfies $\|\mathbf{\varepsilon}_{n+1}\|_{WRMS} \le 1$.  Time steps that satisfy this bound are accepted, while those that violate it are rejected and re-attempted.  In both cases, the step size $h$ is adjusted to strive for the next attempt to satisfy the bound.  We also note that the WRMS norm used for temporal adaptivity is also used in the difference quotient approximation \eqref{eq:difference_quotient} in the next section.

\section{Dominant Eigenvalue Estimation}
\label{sec:dominant-eigenvalue}

In our STS time-integration framework, monitoring the dominant eigenvalue of the Jacobian
operator is essential for computing the required number of stages to maintain numerical stability.
The SUNDIALS library provides a dominant-eigenvalue estimation module, SUNDomEigEst \cite{roberts2025new}, that supports both the power and Arnoldi iterations \cite{trefethen1997numerical, golub2013matrix}. Here, we describe the power iteration approach
used in this work, that is effective for large matrices where the dominant eigenvalue is
real‐valued and has algebraic multiplicity one.

Given a matrix $J$, the power iteration starts with a nonzero initial random vector $\mathbf{v}_0$ and
repeatedly applies
\begin{equation}
\mathbf{v}_{k+1} = \frac{J \mathbf{v}_k}{\|J \mathbf{v}_k\|},
\label{eq:pi-update}
\end{equation}
where $\|\cdot\|$ denotes the Euclidean norm. Over successive iterations, $\mathbf{v}_k$ converges to the
eigenvector associated with the dominant eigenvalue of $J$. At each iteration, the corresponding
eigenvalue is estimated using the Rayleigh quotient
\begin{equation}
\lambda_k = \frac{\mathbf{v}_k^{T} J \mathbf{v}_k}{\mathbf{v}_k^T \mathbf{v}_k}.
\label{eq:pi-rayleigh}
\end{equation}

Importantly, the matrix $J$ need not be formed explicitly; only matrix‐vector products of the form
$J\mathbf{v}$ are required. Thus, power iteration is well suited for PDE discretizations where the Jacobian
is expensive or impractical to assemble, and it introduces minimal memory overhead regardless of the
iteration count.  Since we must estimate the dominant eigenvalues for a potentially varying IVP \eqref{eq:IVP_unsplit}, we approximate these matrix-vector products using SUNDIALS' default difference‐quotient approximation
\begin{equation}
\label{eq:difference_quotient}
J|_{(t,\fvec)} \mathbf{v} \approx \frac{\mathcal{G}(t,\,\fvec + \sigma \mathbf{v}) - \mathcal{G}(t,\,\fvec)}{\sigma},
\end{equation}
where $\sigma = 1/\|\mathbf{v}\|_{WRMS}$ is a small perturbation parameter, and the WRMS norm is defined in either equation \eqref{eq:norm1} or \eqref{eq:norm2}. This enables dominant‐eigenvalue estimation using only evaluations of the right‐hand side function $\mathcal{G}$.

The iterations continue until a relative eigenvalue convergence criterion is satisfied, e.g.,
\begin{equation}
\frac{|\lambda_k - \lambda_{k-1}|}{|\lambda_k|} < \tau,
\label{eq:pi-tolerance}
\end{equation}
for a prescribed tolerance $\tau$. When the dominant eigenvalue is strictly greater in magnitude
than all others, convergence to that eigenvalue is guaranteed, with a linear convergence rate given by the ratio of the magnitudes of the two largest eigenvalues.  Notably, our results in the next section show that at most five iterations are required to converge to tolerance $\tau = 0.1$, amounting to negligible additional computational overhead relative to the timestep solve.

The resulting estimates of the dominant eigenvalue are subsequently multiplied by an ``eigensafety factor,'' $q_{\lambda}$, to define an effective dominant eigenvalue ($\lambda_{eff}=q_{\lambda}\,\lambda_{approx}$) that aims to account for inaccuracies arising from convergence tolerances and difference-quotient approximations -- for $0 < \tau \ll 1$, the eigensafety factor should be at least
\begin{equation}\label{eq:eigsafetybound}
q_{\lambda}>\frac{1}{1-\tau}.
\end{equation}
We investigate this expectation numerically in the next section.

\section{Diffusion Tests}
\label{sec:diffusion-tests}

We test the above methods on the diffusion equation~\eqref{eq:diffusion} over the domain $(t,v_{\|},x)\in[0,t_f]\times[-\pi,\pi]^2$, under periodic boundary conditions, and
where the spatially-dependent diffusion coefficients are given by the function
\begin{equation}
    \label{eq:diffusion_coefficients}
    \begin{split}
       D_{v_{\|}}(v_{\|}) = \nu_{v_{\|}} \left( 1 + 0.99\sin(v_{\|})\right),
    \end{split}
\end{equation}
We note that in gyrokinetic simulations, the diffusion coefficient $\nu_{v_{\|}}$ varies based on the spatial coordinate $x$, thus in the following tests we use $\nu_{v_{\|}} = \left\{0.1, 1, 10\right\}$ to span the range of expected values. The initial condition is given by
\begin{equation}
  \label{eq:diffusion_initial}
  f_0(v_{\|},x) = \frac{1+0.3\sin(2v_{\|})}{\sqrt{5.5\pi}}\;e^{-v_{\|}^2/5.5}.
\end{equation}

To set baseline comparisons between the implicit and explicit time integration methods above, we first test each of the DIRK, SSP and STS methods on a version of this problem that uses a second-order centered finite difference spatial discretization.  As we will see below, these tests indicate that for diffusion equations in this parameter regime, the explicit SSP and STS methods are competitive (and frequently outperform) the DIRK methods.  Therefore, in Section \ref{sec:dg-results} we transition to the DG spatial discretization, and focus on only the SSP and STS methods since they do not require construction of Jacobian entries.

\subsection{STS, SSP and DIRK Comparison for the Finite Difference Discretization}
\label{sec:fd-results}

We first test the DIRK, SSP, and STS methods on the diffusion equation \eqref{eq:diffusion} with diffusion coefficient \eqref{eq:diffusion_coefficients} and initial condition \eqref{eq:diffusion_initial}, using a second-order centered finite difference spatial discretization.  We discretize the spatial domain $[-\pi,\pi]^2$ using a uniform mesh with $N_{v_\|} \times N_x$ grid points, where $N_{v_\|} = N_x = \{64, 128, 256\}$, where as we refine the grid we parallelized with a grid of $\left\{2\times 2, 4\times 4, 8\times 8\right\}$ MPI ranks, respectively.  The final time is set to $t_f = 1.0$, and we use diffusion coefficient $\nu_{v_\|} = \{0.1,1,10\}$.  We test each time integration method using relative tolerances of $RTOL = 10^{-k}, k=2,\ldots,6$, with an absolute tolerance of $ATOL = 10^{-11}$.  We compare the computational efficiency of each method by plotting the relative temporal error against the runtime -- thus a method with results to the left of other methods for a given accuracy is considered more efficient.  Here, we determine accuracy by measuring the $L^{\infty}$ relative error at a set of 20 evenly-spaced times throughout the simulation, comparing against a reference solution computed using $RTOL = 10^{-12}$.

\begin{figure}[H]
    \centering
    \includegraphics[trim={0 0 110 0}, clip, width=0.9\textwidth]{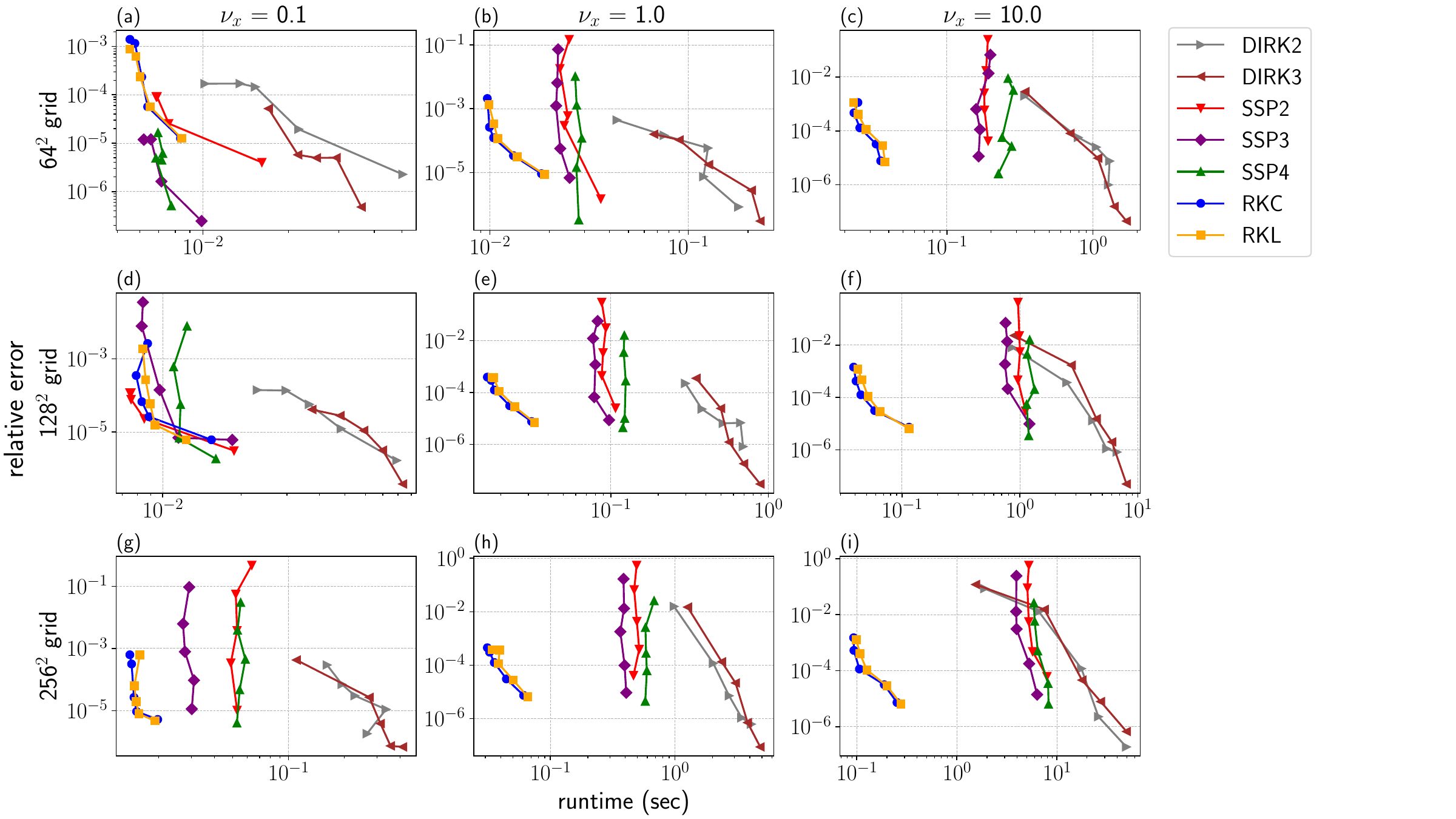}
    \caption{Computational efficiency for each adaptive solver for various tolerances, grids, and diffusion coefficients.}
    \label{fig:fd_results_adaptivity}
\end{figure}

Our efficiency plots for all grids and diffusion constants are shown in Figure \ref{fig:fd_results_adaptivity}.  Recalling that we requested relative errors in the range $[10^{-6}, 10^{-2}]$, we observe that the relative errors obtained from both the SSP and DIRK methods vary considerably based on the test configuration, sometimes resulting in solution errors well above or below the requested relative tolerance ($RTOL$) values, whereas the STS methods consistently achieve solution errors within a factor of 10 from the specified $RTOL$.  Second, we see that the relative performance of the SSP methods depend significantly on the diffusion strength -- for $\nu_{v_\|}=0.1$ the SSP methods are among the most efficient (particularly for coarse grids), but as $\nu_{v_\|}$ increases the SSP performance degrades.  Lastly, we see that in all test configurations considered here, the DIRK methods are considerably less efficient than STS.  Thus, due to their added complication of requiring Jacobian information, we do not include DIRK methods in the following tests using the DG spatial discretization.

\subsection{STS and SSP Comparison for the DG Discretization}
\label{sec:dg-results}

In this section, we present our implementation of the STS and SSP methods for the diffusion equation using a DG discretization with $N_{v_{\parallel}}=120$ and $N_x=20$ cells and a piecewise linear basis, and compare their computational performance. As discussed in Section \ref{sec:STSmethods}, STS methods predetermine the number of stages -- and hence the stability region -- prior to taking each time step, based on a prescribed estimate of the dominant eigenvalue ($\lambda_{\text{user}}$ or $\lambda_{\text{approx}}$). We investigate the impact of these two choices on the numerical solution below. In addition, we examine the performance of our two WRMS norm variants (equations \eqref{eq:norm1} and \eqref{eq:norm2}) for temporal error estimation.

We first consider how to identify a reasonable eigensafety value ($q_{\lambda}$), since too small a value could lead to numerical instabilities (and thus rejected time steps), and too large a value would lead to excessive computational effort.  We note that in general, user-provided formulas for the dominant eigenvalues $\lambda_{user}$ over-estimate the true value, particularly in comparison with direct measurements provided by SUNDIALS the internal dominant eigenvalue estimation, $\lambda_{approx}$. For the latter, we fix the tolerance for estimating $\lambda_{approx}$ to $\tau = 0.1$.  Thus according to our estimate \eqref{eq:eigsafetybound}, a safety factor of at least $1/0.9$ should be used to guarantee numerical stability. 

\begin{figure}[H]
    \centering
    \begin{subfigure}[b]{0.355\textwidth}
        \includegraphics[trim={0 0 0 0}, clip, width=\textwidth]{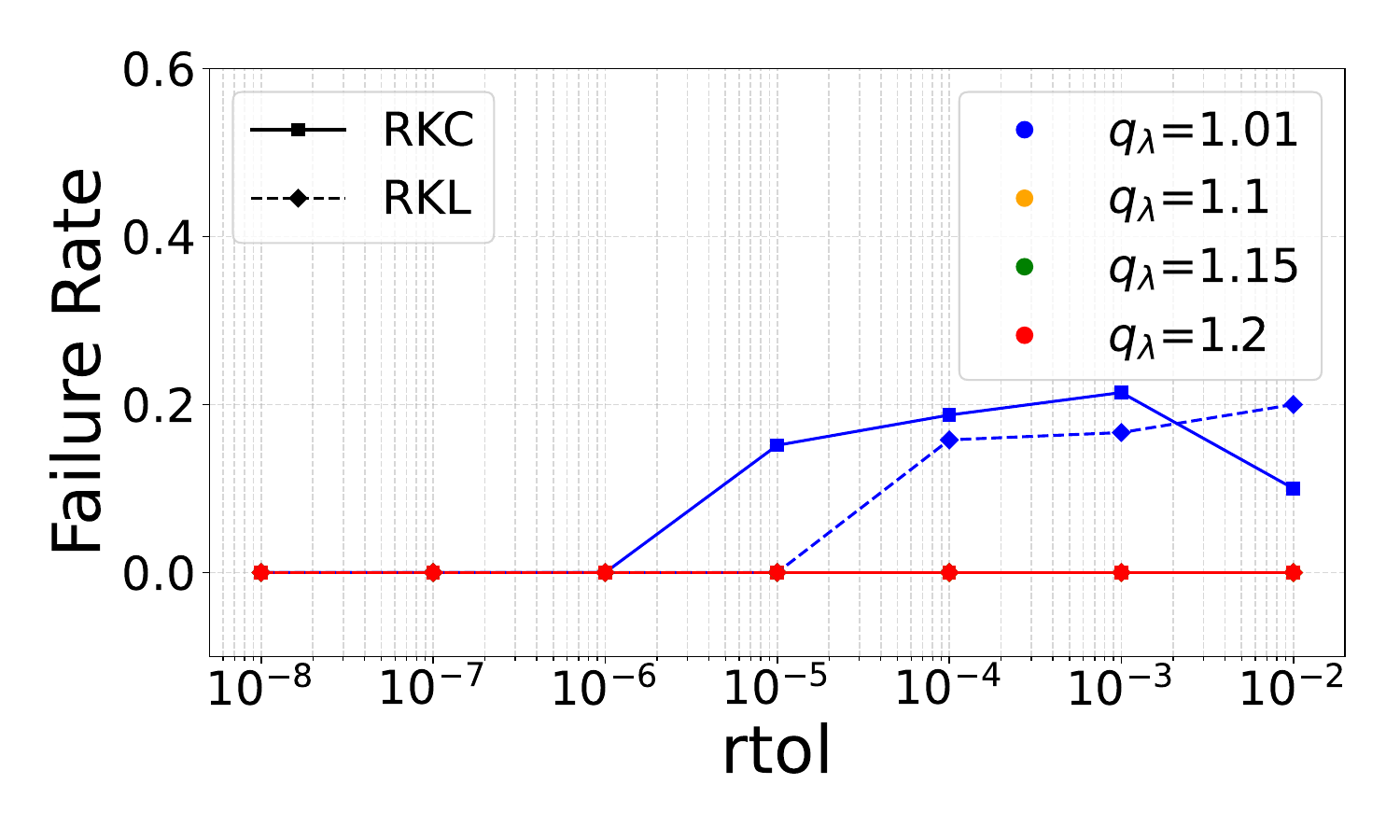}
        \caption{$\nu_{v_\|}=0.1$}
        \label{fig:eigensafety_rtol_vs_FR1}
    \end{subfigure}
    \hfill
    \begin{subfigure}[b]{0.305\textwidth}
        \includegraphics[trim={100 0 0 0}, clip, width=\textwidth]{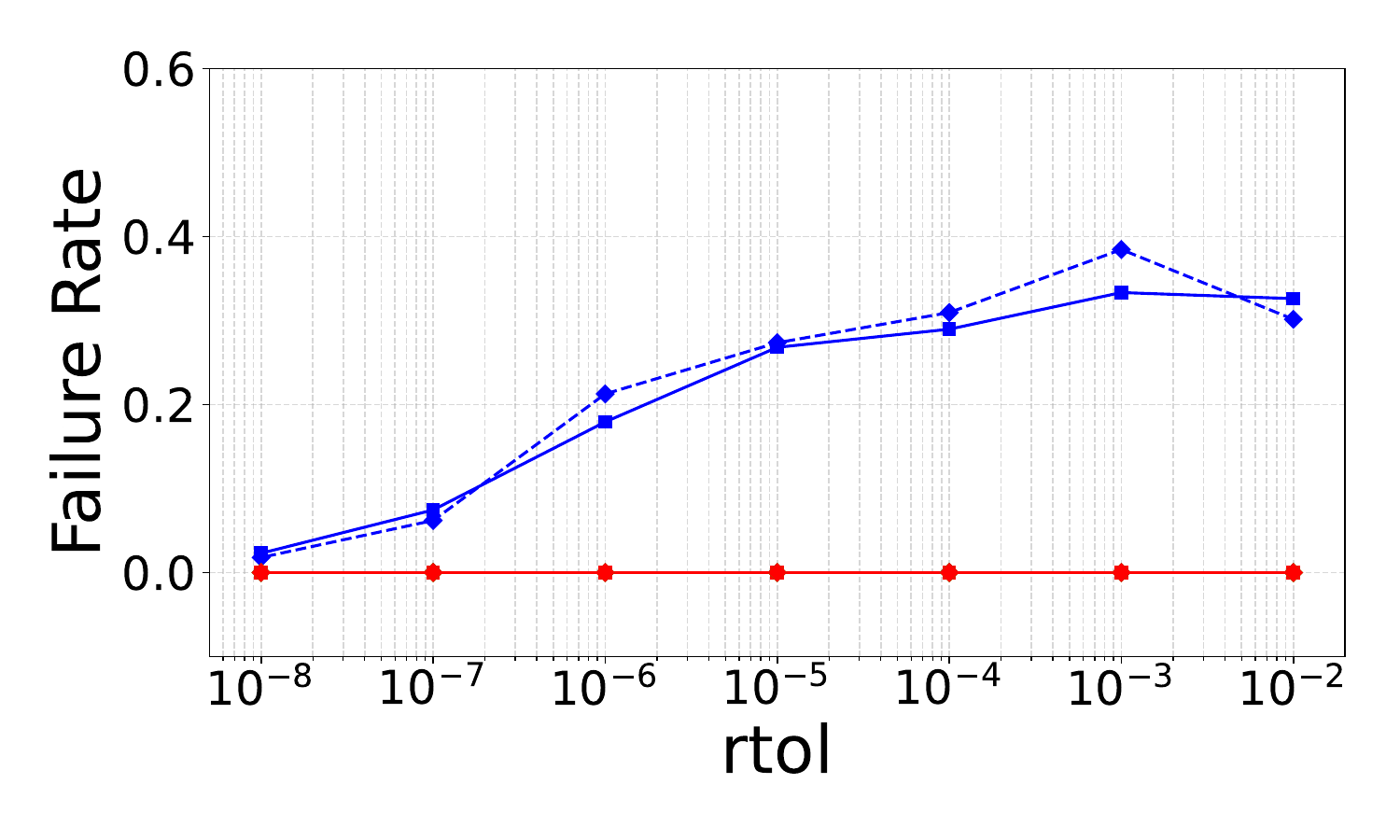}
        \caption{$\nu_{v_\|}=1.0$}
        \label{fig:eigensafety_rtol_vs_FR2}
    \end{subfigure}
    \hfill
    \begin{subfigure}[b]{0.305\textwidth}
        \includegraphics[trim={100 0 0 0}, clip, width=\textwidth]{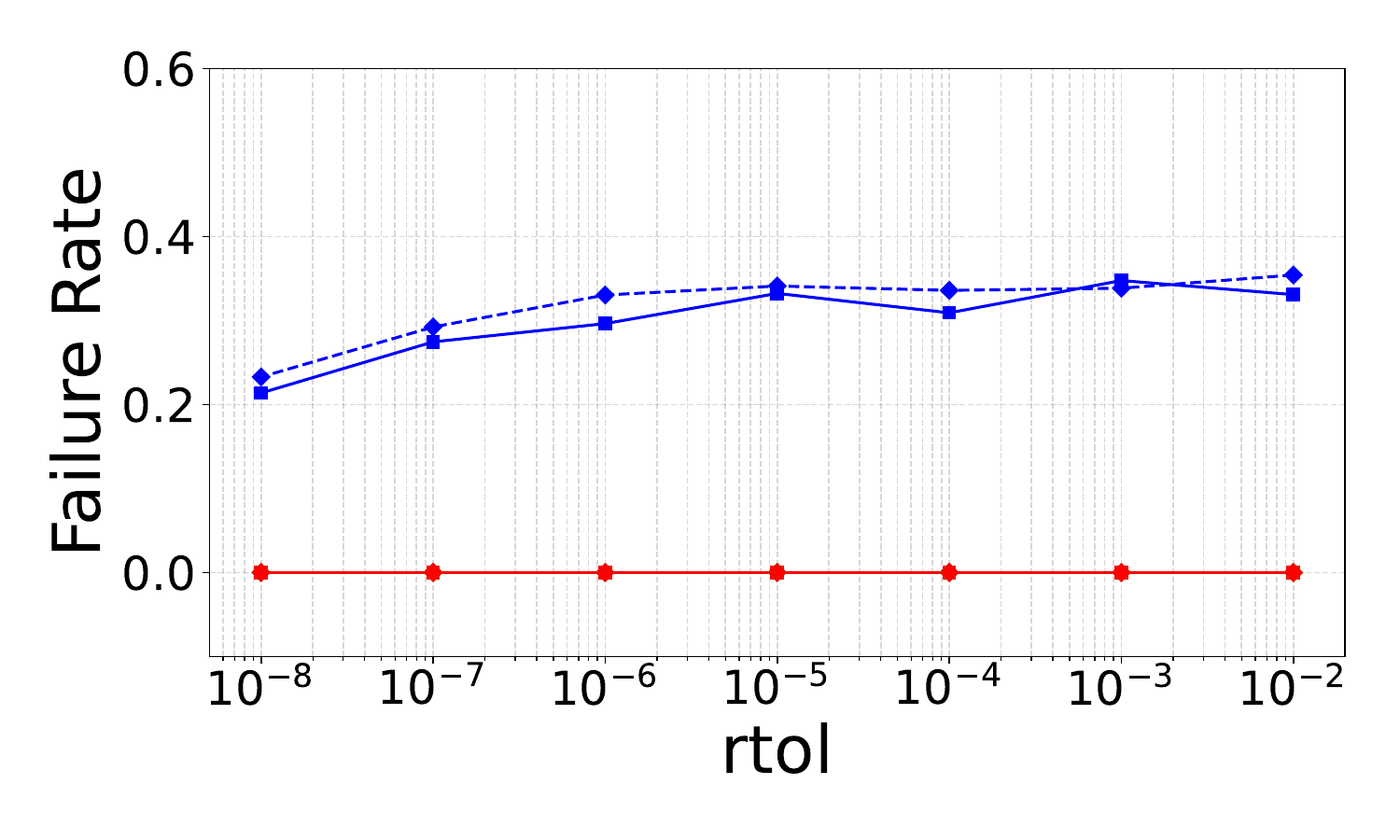}
        \caption{$\nu_{v_\|}=10.0$}
        \label{fig:eigensafety_rtol_vs_FR3}
    \end{subfigure}

    \caption{Sensitivity to dominant eigenvalue safety factor ($q_{\lambda}$)}
    \label{fig:eigensafety_rtol_vs_FR}
\end{figure}

Figure \ref{fig:eigensafety_rtol_vs_FR} illustrates the step failure rates for each relative tolerance values $(\text{rtol} = 10^{-i},\, i=2,\dots,8,$ therein and the rest of the computations). It suggest that safety factors below $1.1$ lead to failure rates up to 20\%, which is in line with our expected bound from equation \eqref{eq:eigsafetybound}. Interestingly, however, the safety factor $q_{\lambda}=1.1$ results in zero failures, even though it is slightly smaller than our estimated minimum of $1/0.9$.  We expect that this value is sufficient because when calculating the minimum number of STS stages, $s$ must be rounded up to the nearest integer.  Thus in the remaining tests we use $q_{\lambda}=1.1$, but caution that for other applications a more conservative safety factor may be prudent.

\begin{figure}[H]
    \centering
    \begin{subfigure}[b]{0.371\textwidth}
        \includegraphics[trim={0 0 0 0}, clip, width=\textwidth]{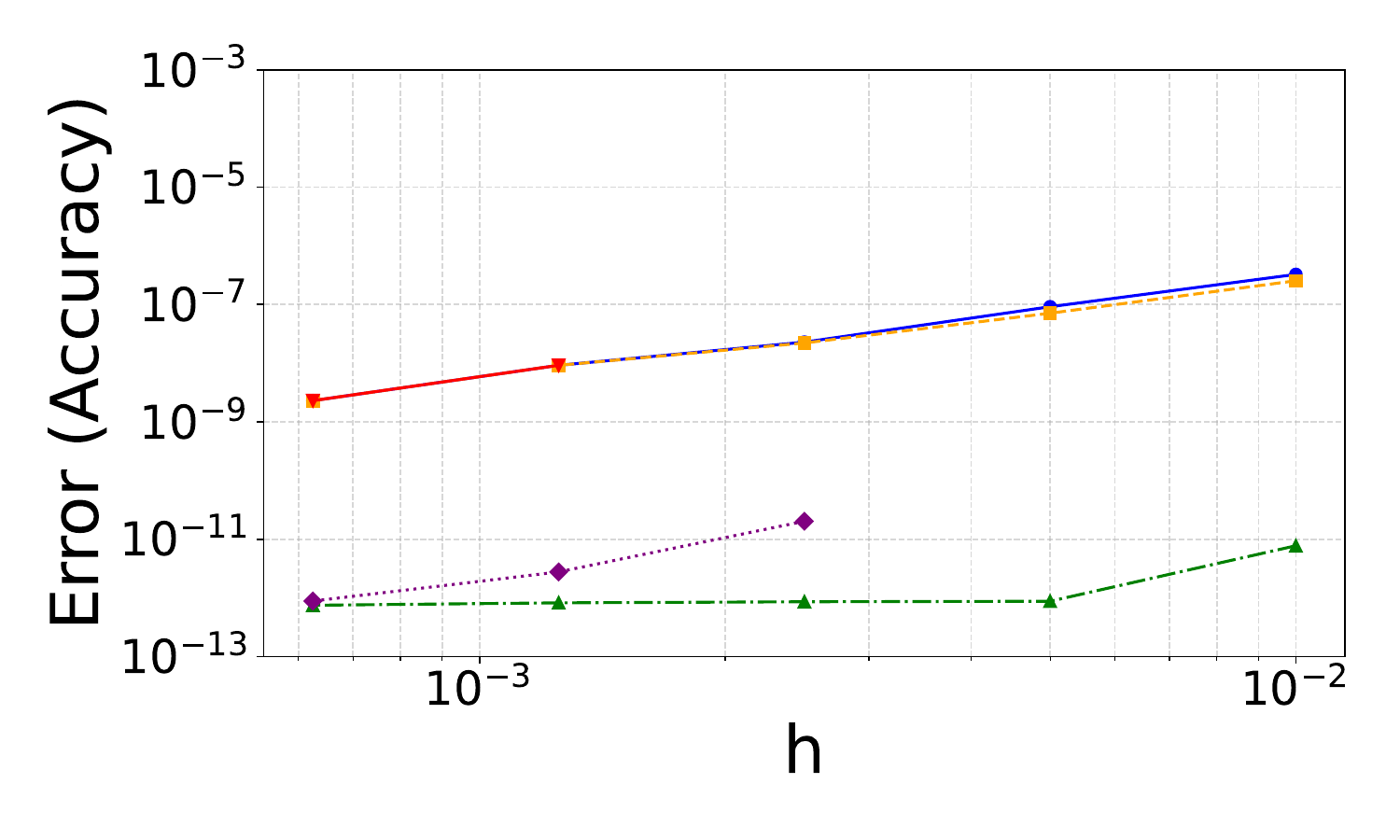}
        \caption{$\nu_{v_\|}=0.1$}
        \label{fig:SSP_fails1}
    \end{subfigure}
    \hfill
    \begin{subfigure}[b]{0.305\textwidth}
        \includegraphics[trim={128 0 0 0}, clip, width=\textwidth]{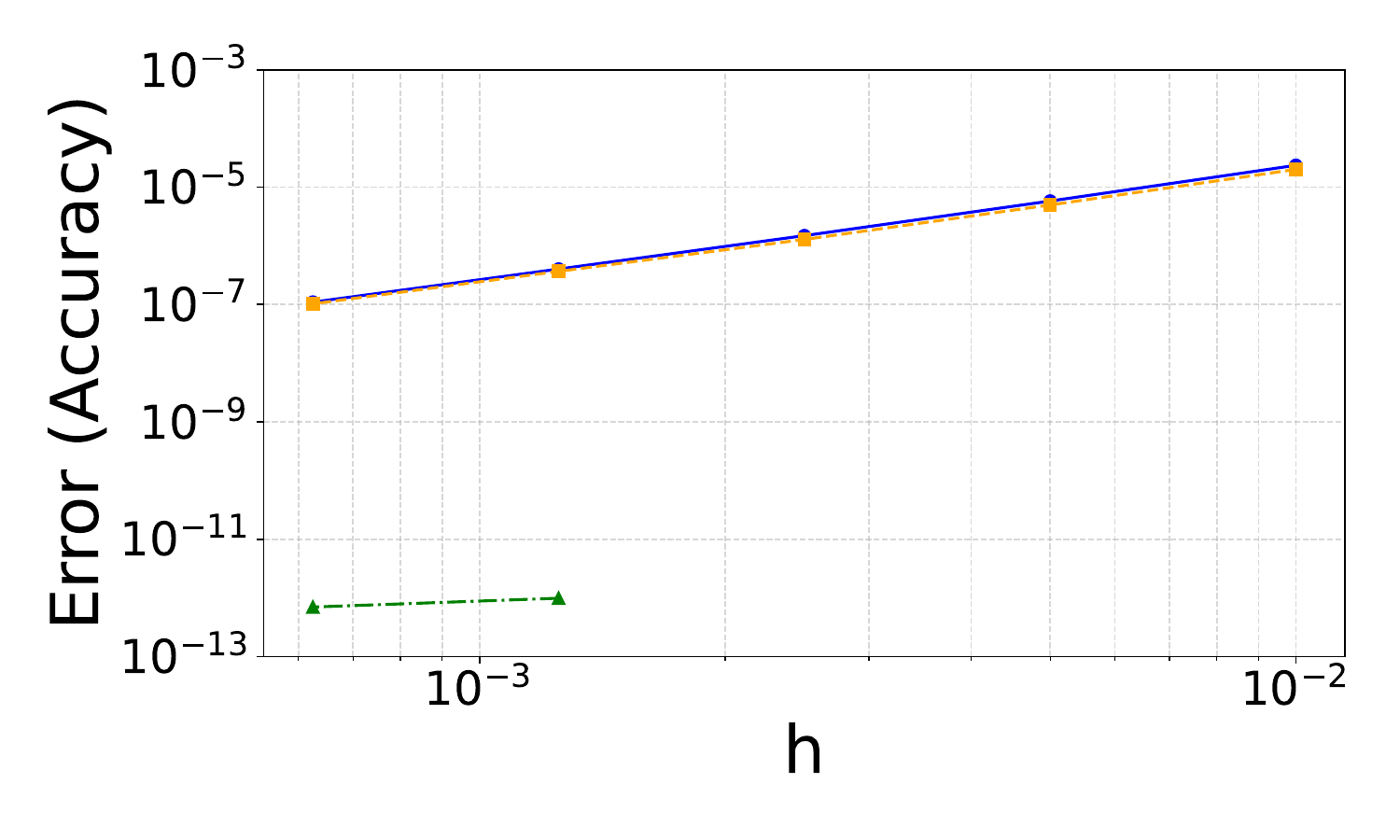}
        \caption{$\nu_{v_\|}=1.0$}
        \label{fig:SSP_fails2}
    \end{subfigure}
    \hfill
    \begin{subfigure}[b]{0.305\textwidth}
        \includegraphics[trim={128 0 0 0}, clip, width=\textwidth]{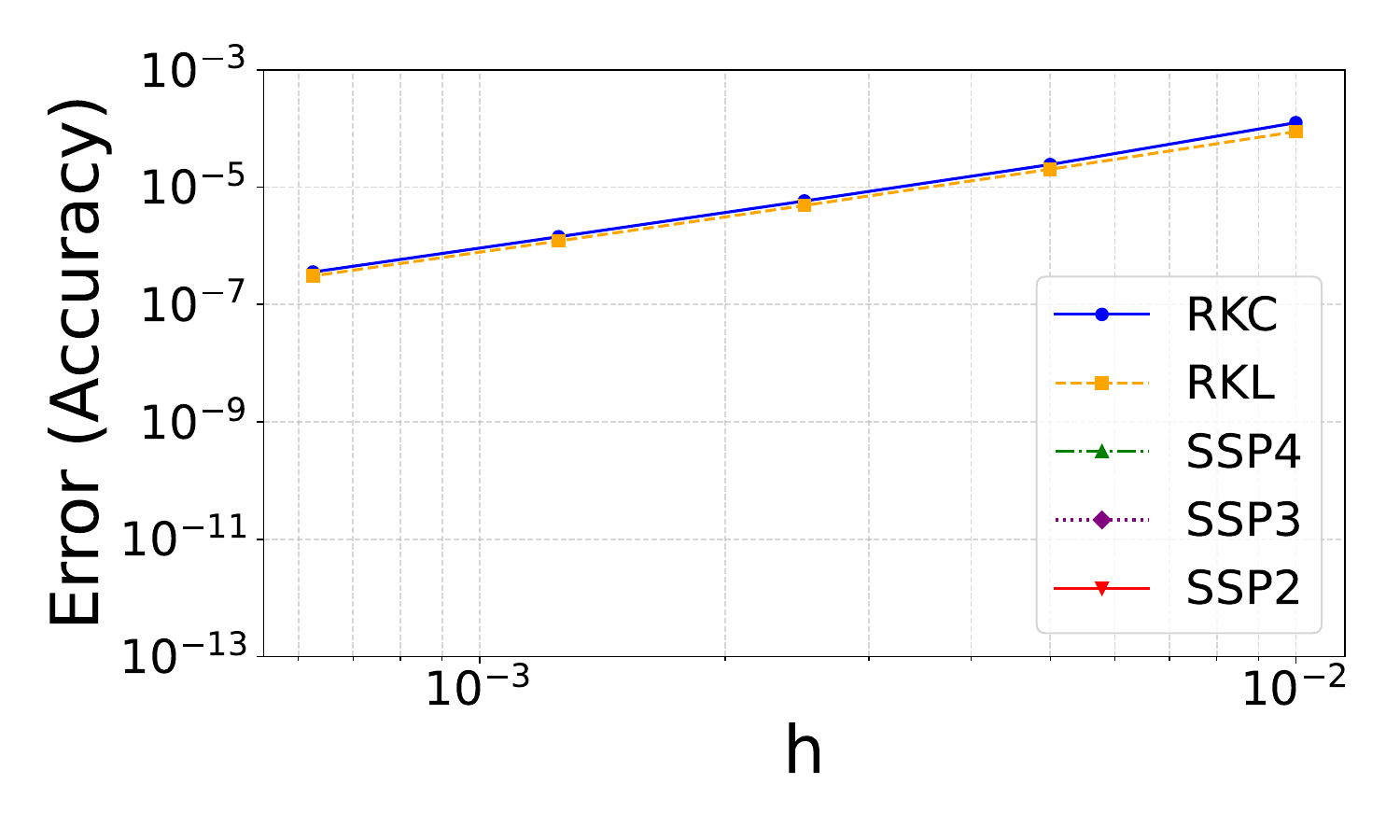}
        \caption{$\nu_{v_\|}=10.0$}
        \label{fig:SSP_fails3}
    \end{subfigure}

    \caption{Error vs fixed step size ($h$) for varying diffusion coefficient $\nu_{v_\|}$}
    \label{fig:SSP_fails}
\end{figure}

Next, we investigate the stability properties of STS and SSP methods. Both methods are explicit and hence subject to stability constraints. In Figure~\ref{fig:SSP_fails} we plot the solution error versus fixed step size $h$ (chosen as 0.01 for the first refinement level and then successively halved 4 times) for each method and diffusion coefficient $\nu_{v\|}$, omitting results that blew up due to numerical instability.  This illustrates that the variable stability of the STS methods allows it to take time steps that far exceed the largest stable steps for the SSP methods, while continuing to provide their predicted second order convergence.  We note that both SSP2 and SSP3 were only stable for the smallest steps in Figure \ref{fig:SSP_fails1}, and although SSP4 produced results for all steps in Figure \ref{fig:SSP_fails1}, it also becomes unstable as the diffusion coefficient $\nu_{v_\|}$ increases in Figures \ref{fig:SSP_fails2} and \ref{fig:SSP_fails3}.

Our next tests compare the two competing WRMS norms from equations \eqref{eq:norm1} and \eqref{eq:norm2}, that correspond with SUNDIALS' default component-wise norm that measures relative error in each solution component, and the proposed cell-wise norm instead that measures relative error in each DG element.  As stated in Section \ref{sec:TemporalAdaptivity}, we predict the latter will be beneficial for DG discretizations where the degrees of freedom in an element contribute unequally to the overall solution.

\begin{figure}[H]
    \centering
    \begin{subfigure}[b]{0.371\textwidth}
        \includegraphics[trim={0 0 0 0}, clip, width=\textwidth]{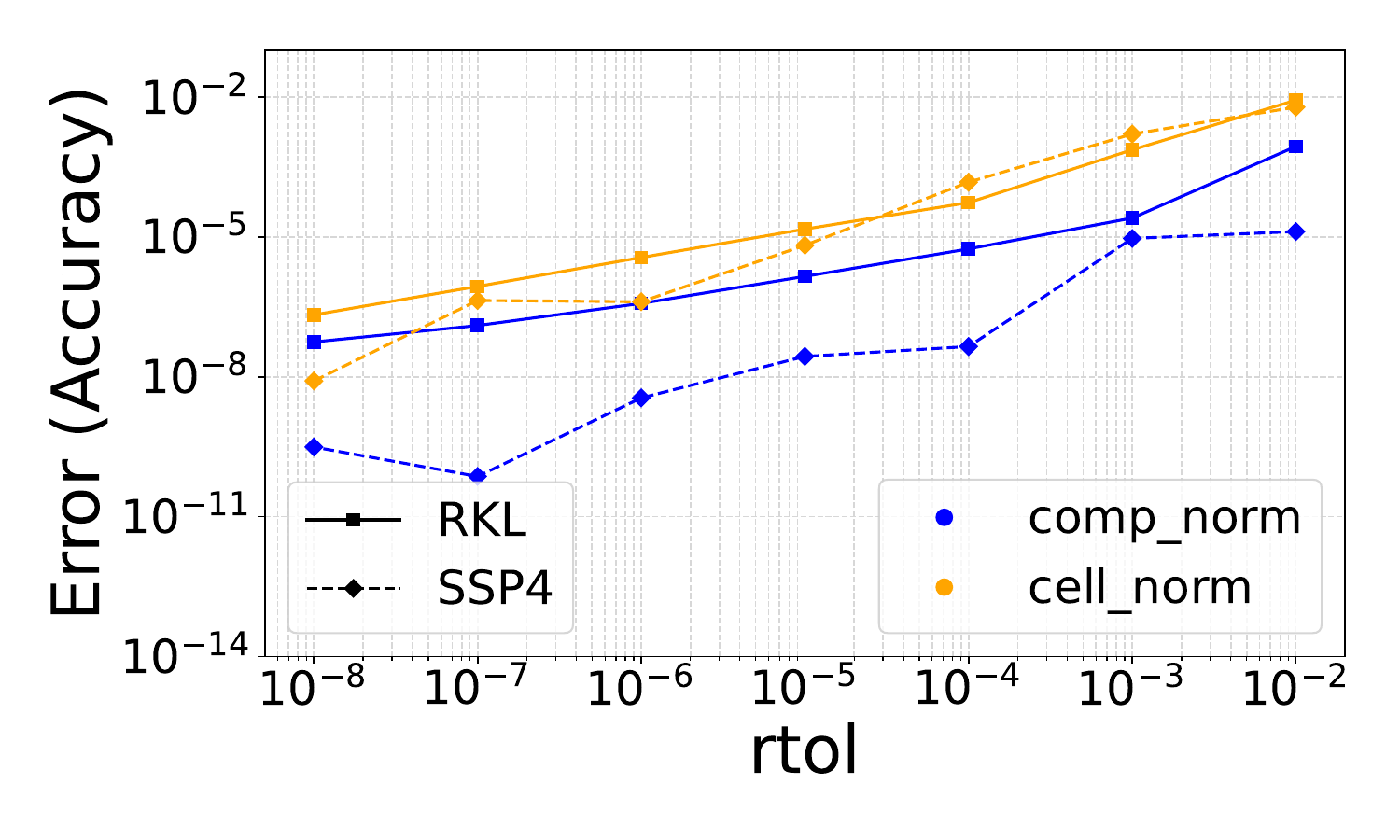}
        \caption{$\nu_{v_\|}=0.1$}
        \label{fig:error_vs_rtol1}
    \end{subfigure}
    \hfill
    \begin{subfigure}[b]{0.305\textwidth}
        \includegraphics[trim={128 0 0 0}, clip, width=\textwidth]{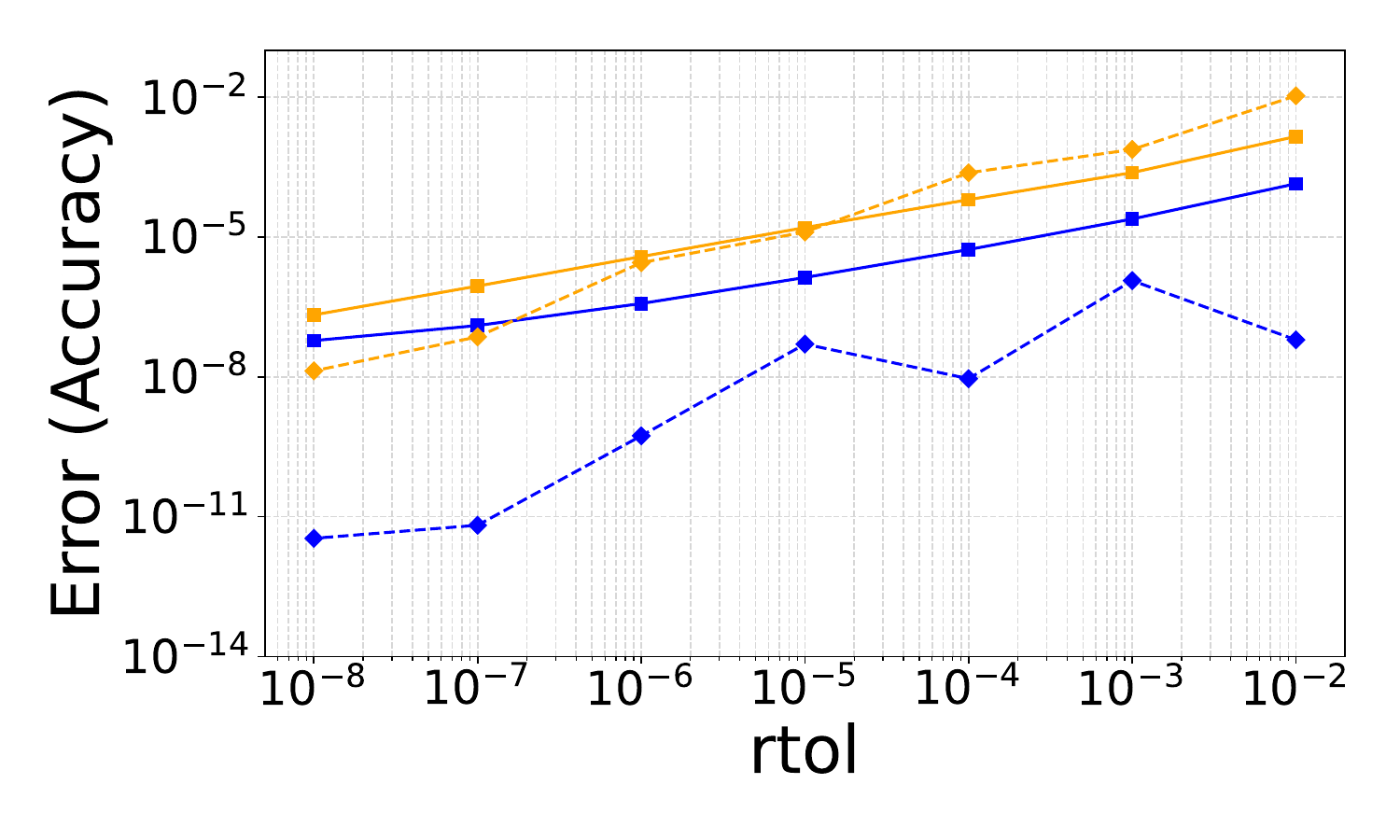}
        \caption{$\nu_{v_\|}=1.0$}
        \label{fig:error_vs_rtol2}
    \end{subfigure}
    \hfill
    \begin{subfigure}[b]{0.305\textwidth}
        \includegraphics[trim={128 0 0 0}, clip, width=\textwidth]{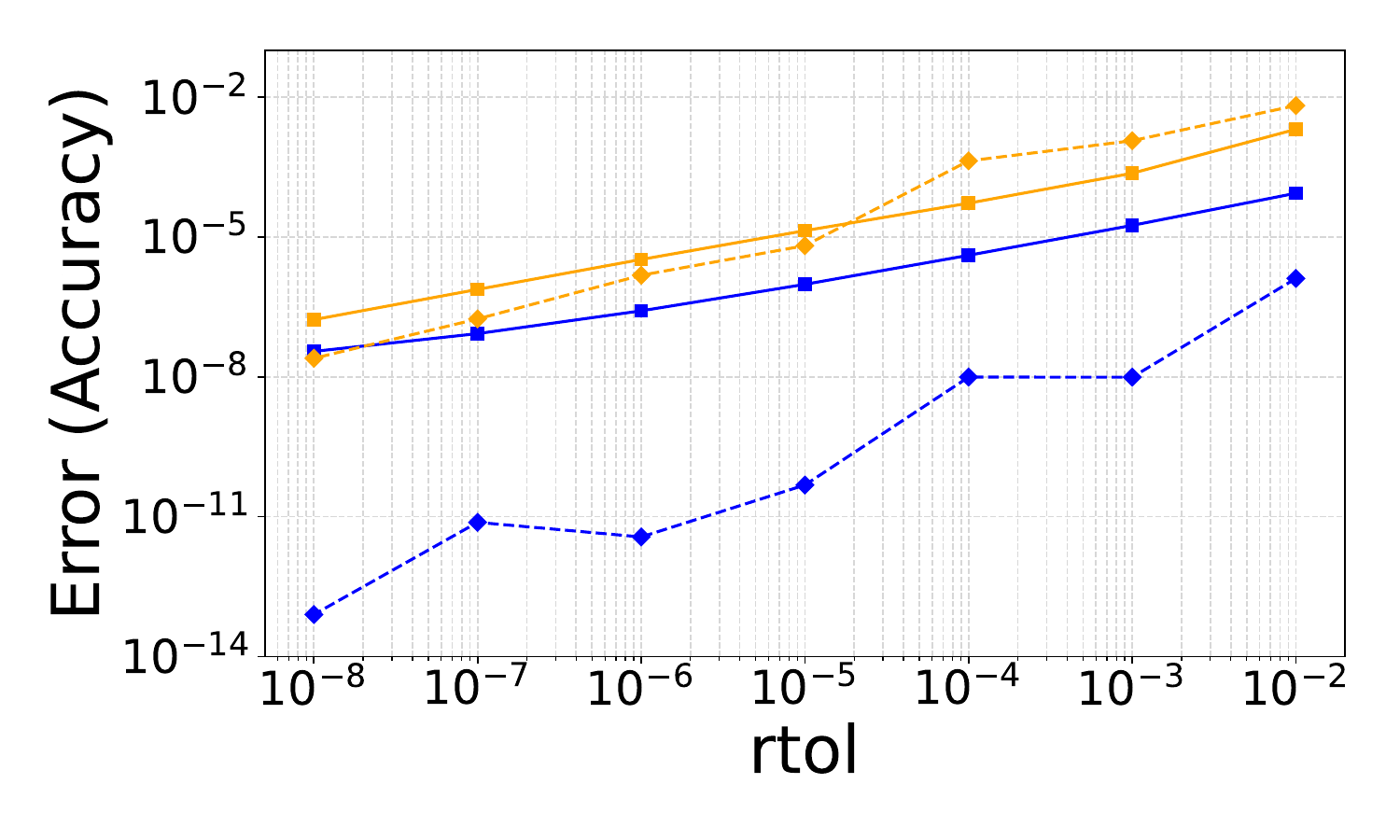}
        \caption{$\nu_{v_\|}=10.0$}
        \label{fig:error_vs_rtol3}
    \end{subfigure}

    \caption{Error vs relative tolerance (rtol) for the component-wise and cell-wise norms, equations \eqref{eq:norm1} and \eqref{eq:norm2} respectively, when using RKL and SSP4.}
    \label{fig:error_vs_rtol}
\end{figure}

Recall that temporal adaptivity strategies select time steps such that the local error in each step satisfies $\|\mathbf{\varepsilon}_{n+1}\|_{WRMS} \le 1$, where the norm incorporates user-requested relative and absolute tolerances.  However, due to error accumulation between steps, it is expected that the overall solution error can exceed these tolerances.  Figure \ref{fig:error_vs_rtol} displays the solution error measured in the max-in-time, max-in-space norm,
\begin{equation}
\label{eq:maxmaxnorm}
\max_{1 \le n \le N} \|\fvec_{n} - \fvec_{\text{true}}(t_n)\|_{\ell^\infty},
\end{equation}
versus the relative tolerance for RKL and SSP4 (the RKC, SSP2 and SSP3 methods exhibit similar behavior and are therefore not shown to avoid overcrowding the figures).

Figure \ref{fig:error_vs_rtol} indicates that when using the component-wise norm, SSP4 produces solutions that are much more accurate than requested.  While low error is not problematic in itself, this additional accuracy can come with additional costs, including (i) increased computation time, (ii) additional step rejections (steps that do not meet the requested error tolerance).

\begin{figure}[H]
    \centering
    \begin{subfigure}[b]{0.371\textwidth}
        \includegraphics[trim={0 0 0 0}, clip, width=\textwidth]{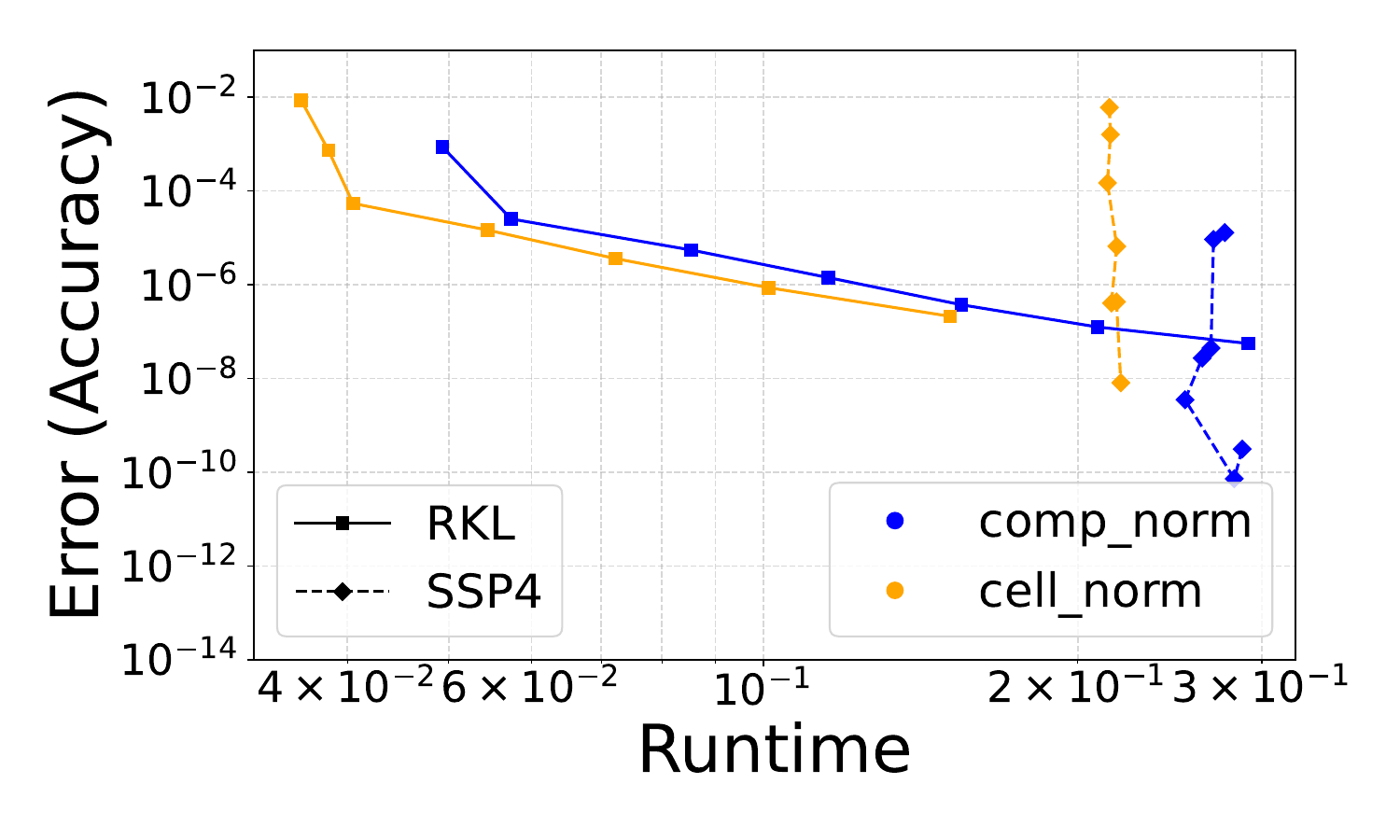}
        \caption{$\nu_{v_\|}=0.1$}
        \label{fig:error_vs_runtime1}
    \end{subfigure}
    \hfill
    \begin{subfigure}[b]{0.305\textwidth}
        \includegraphics[trim={128 0 0 0}, clip, width=\textwidth]{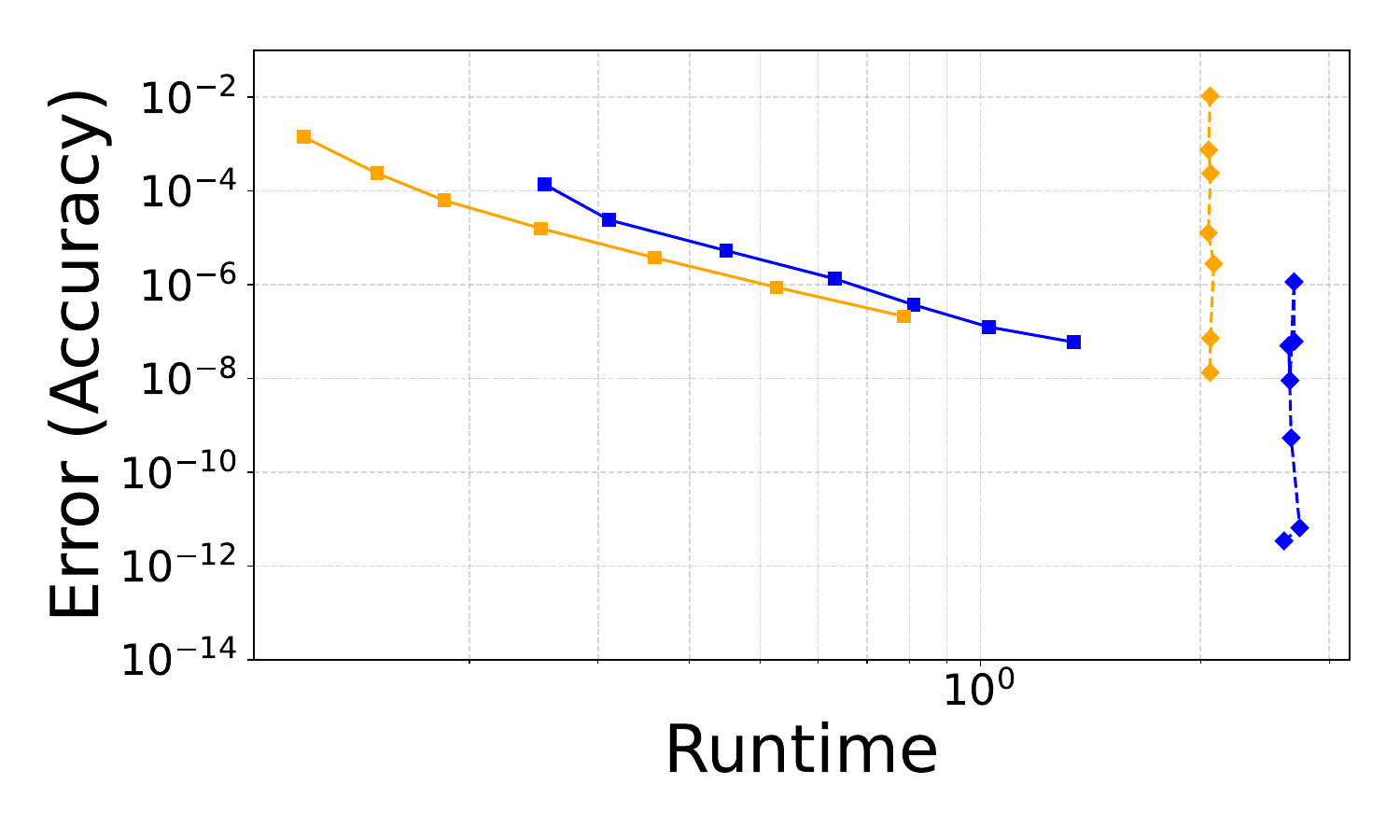}
        \caption{$\nu_{v_\|}=1.0$}
        \label{fig:error_vs_runtime2}
    \end{subfigure}
    \hfill
    \begin{subfigure}[b]{0.305\textwidth}
        \includegraphics[trim={128 0 0 0}, clip, width=\textwidth]{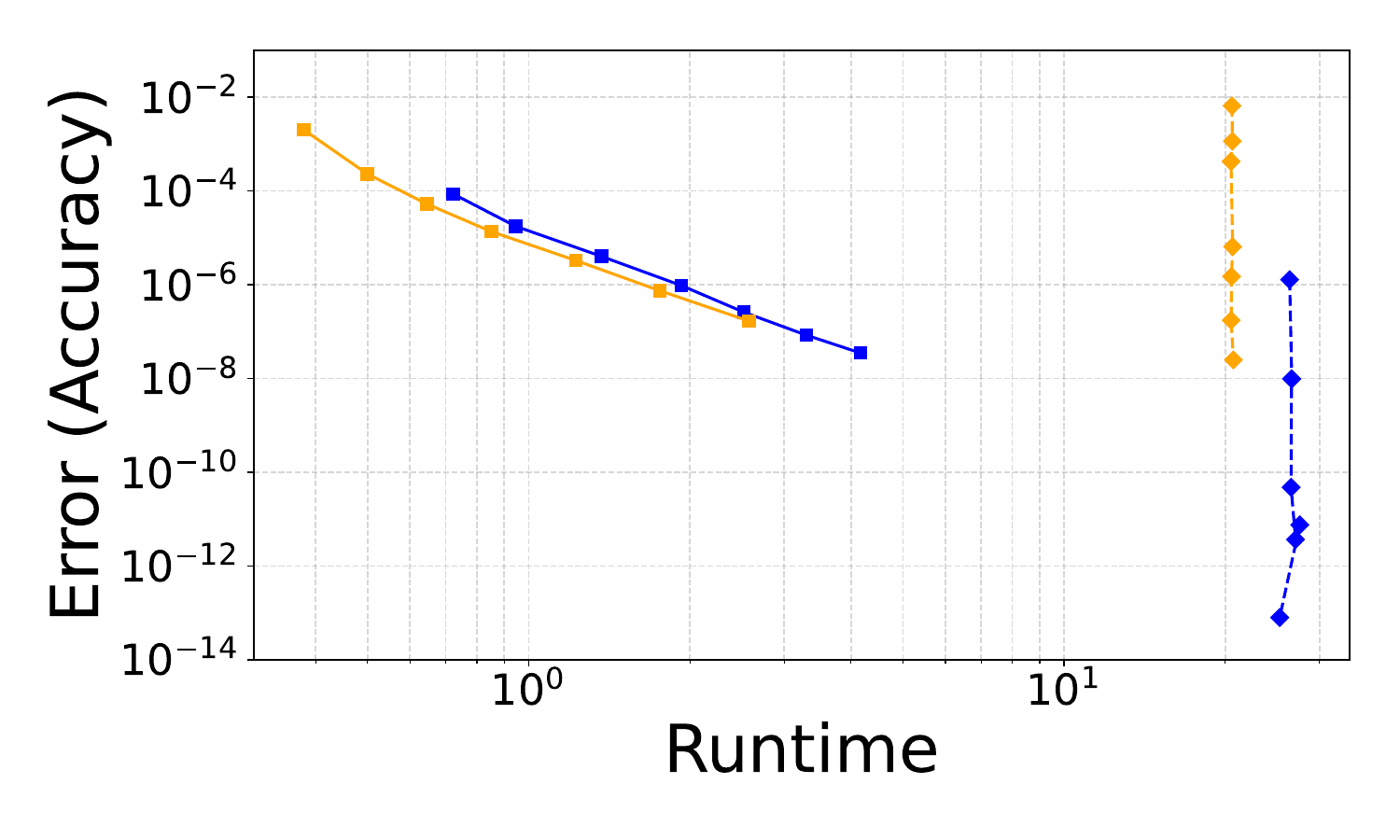}
        \caption{$\nu_{v_\|}=10.0$}
        \label{fig:error_vs_runtime3}
    \end{subfigure}

    \caption{Error vs computation times for the component-wise and cell-wise norms, equations \eqref{eq:norm1} and \eqref{eq:norm2} respectively, when using RKL or SSP4.}
    \label{fig:error_vs_runtime}
\end{figure}

To assess these increased costs, Figure \ref{fig:error_vs_runtime} plots the computation times to obtain corresponding accuracies.  As in Figure \ref{fig:fd_results_adaptivity}, we consider a method to be more efficient if its curve lies to the left of competing methods, in that it requires less runtime to achieve the same solution accuracy. Here we see that for both the RKL and SSP4 methods, the cell-wise norm is demonstrably more efficient than the component-wise norm, and that for nearly all solution accuracies, RKL is more efficient than SSP4.  Additionally, for RKL the computational efficiency significantly drops when the component-wise norm is employed with tight tolerances.

\begin{figure}[H]
    \centering
    \begin{subfigure}[b]{0.355\textwidth}
        \includegraphics[trim={0 0 0 0}, clip, width=\textwidth]{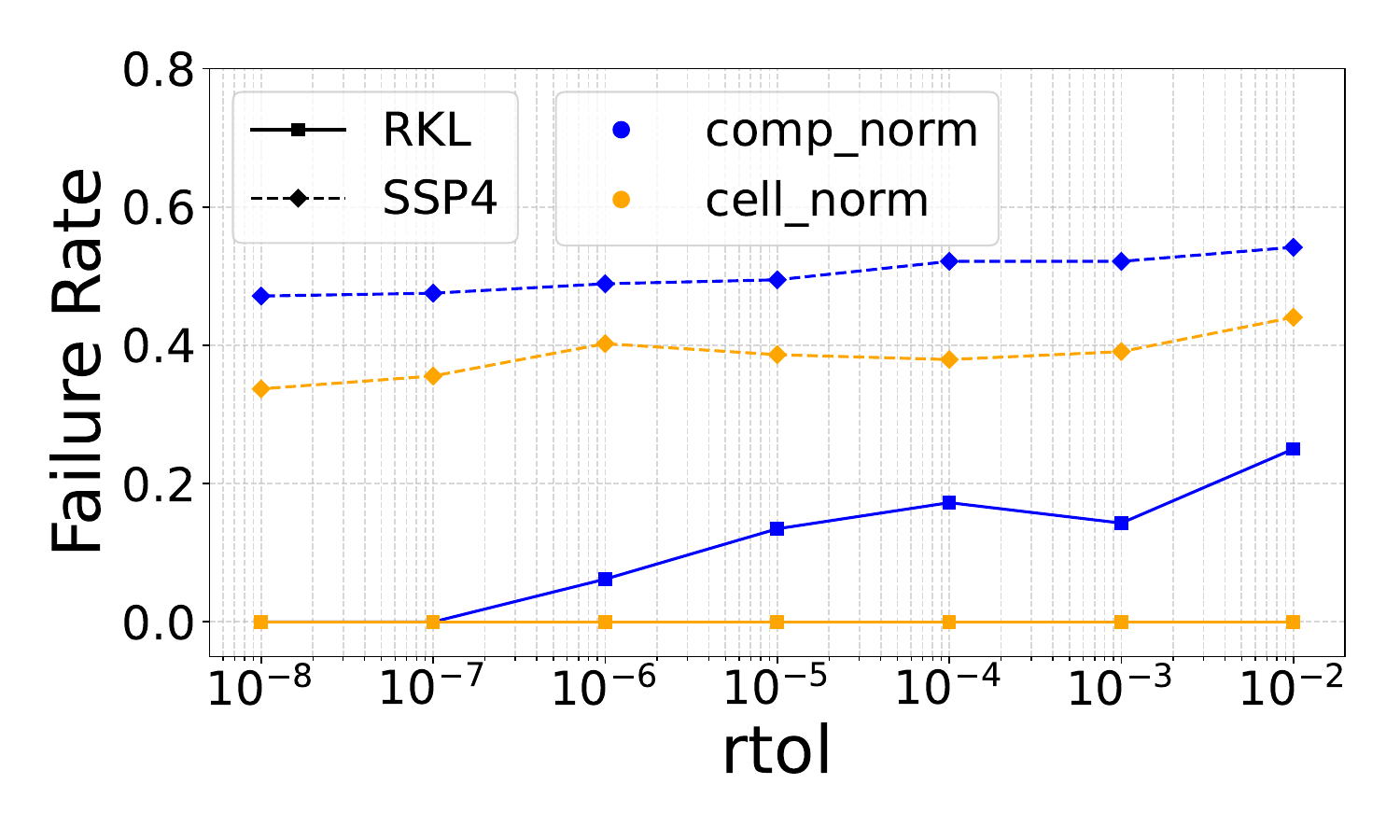}
        \caption{$\nu_{v_\|}=0.1$}
        \label{fig:rtol_vs_FR1}
    \end{subfigure}
    \hfill
    \begin{subfigure}[b]{0.305\textwidth}
        \includegraphics[trim={100 0 0 0}, clip, width=\textwidth]{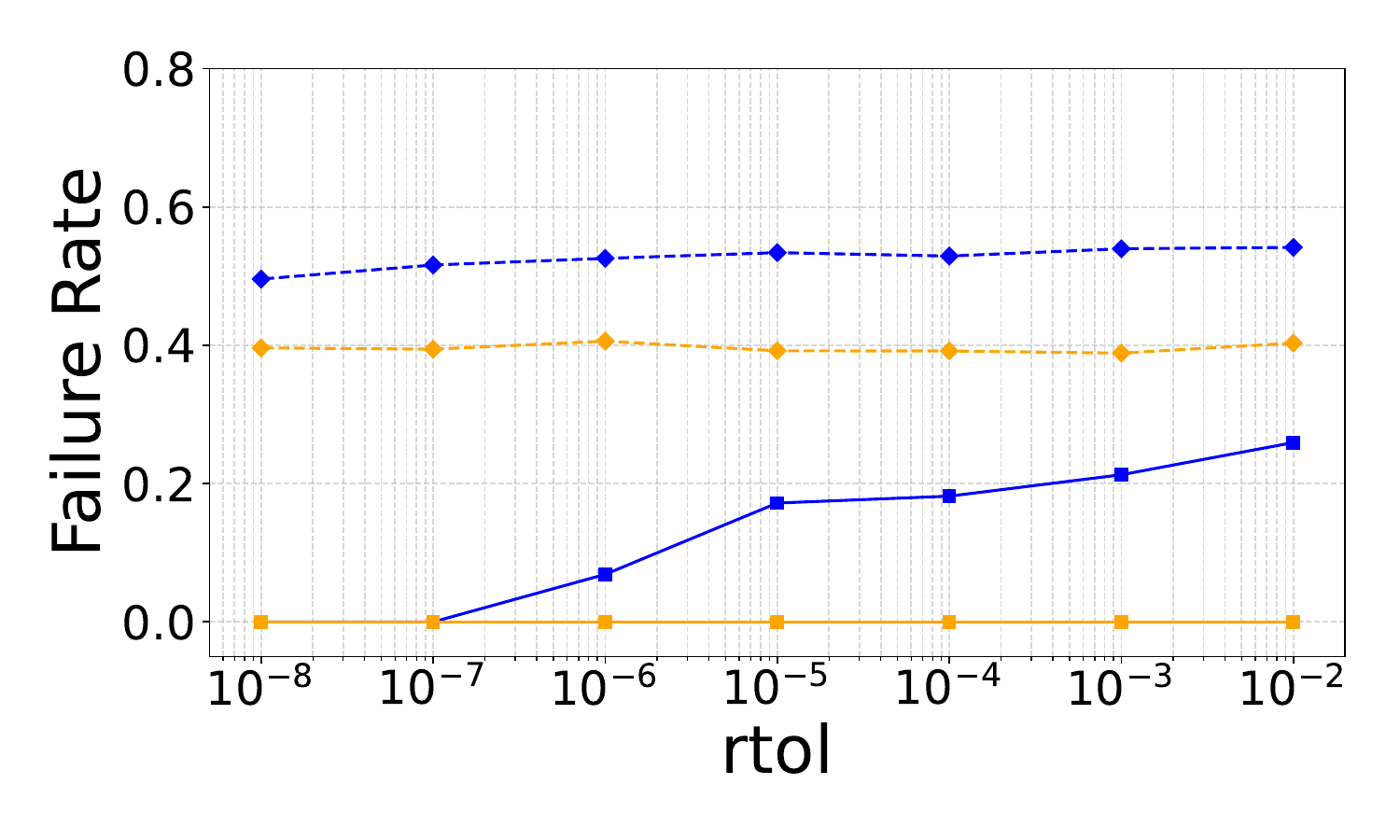}
        \caption{$\nu_{v_\|}=1.0$}
        \label{fig:rtol_vs_FR2}
    \end{subfigure}
    \hfill
    \begin{subfigure}[b]{0.305\textwidth}
        \includegraphics[trim={100 0 0 0}, clip, width=\textwidth]{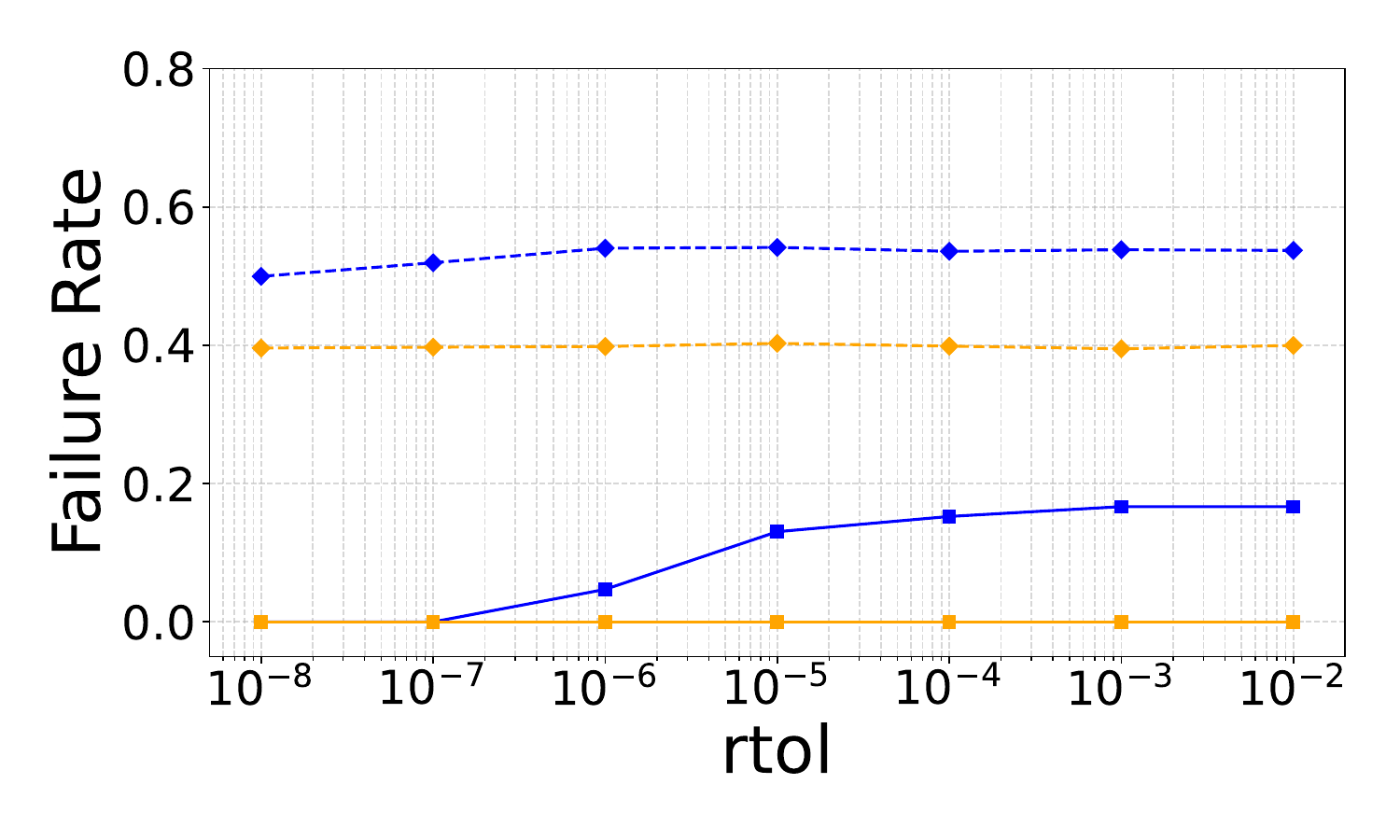}
        \caption{$\nu_{v_\|}=10.0$}
        \label{fig:rtol_vs_FR3}
    \end{subfigure}

    \caption{Step failure rates vs relative tolerance  for differing norm types}
    \label{fig:rtol_vs_FR}
\end{figure}

To further understand this performance difference, we may consider the percentage of steps that were rejected by each time adaptivity algorithm, since these failed steps contribute directly to computational cost.  Figure \ref{fig:rtol_vs_FR} plots these failure rates versus the relative tolerance, illustrating that the RKL method (similarly for RKC) has no rejected steps when using the cell-wise norm. The component-wise norm with loose tolerances increases the rate of the rejected STS steps up to 26\%. This observation, combined with the observations from Figure \ref{fig:error_vs_rtol}, suggests that the component-wise norm attempts to ``over-solve'' the problem, rejecting steps that would have otherwise have resulted in sufficient accuracy, and in turn resulting in overly-accurate solutions.  For both norms, the SSP method results in many step rejections due to its limited stability region. Nevertheless, the cell-wise norm reduces the failure rate for this method as well.

\begin{figure}[H]
    \centering
    \begin{subfigure}[b]{0.366\textwidth}
        \includegraphics[trim={0 0 0 0}, clip, width=\textwidth]{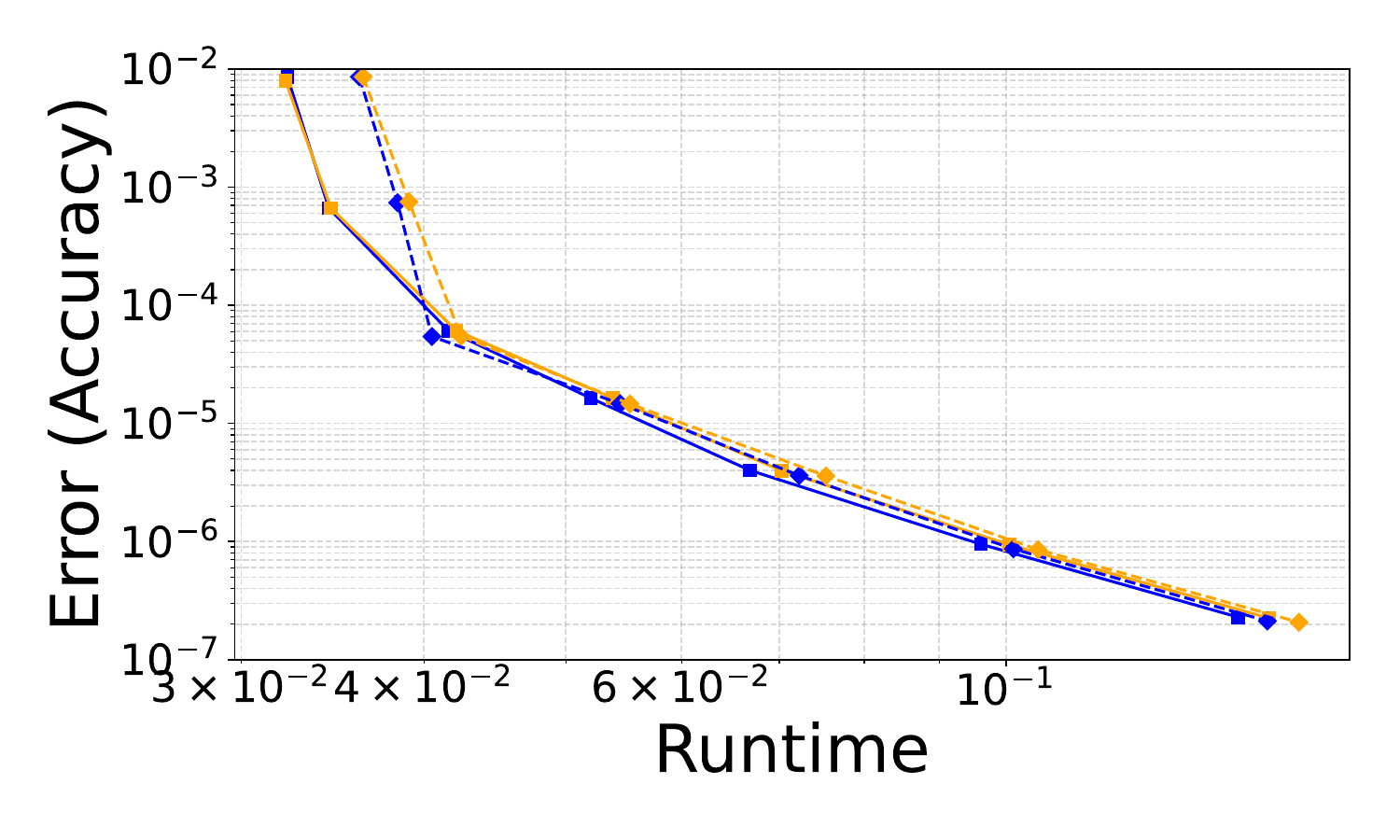}
        \caption{$\nu_{v_\|}=0.1$}
        \label{fig:error_runtime1}
    \end{subfigure}
    \hfill
    \begin{subfigure}[b]{0.305\textwidth}
        \includegraphics[trim={120 0 0 0}, clip, width=\textwidth]{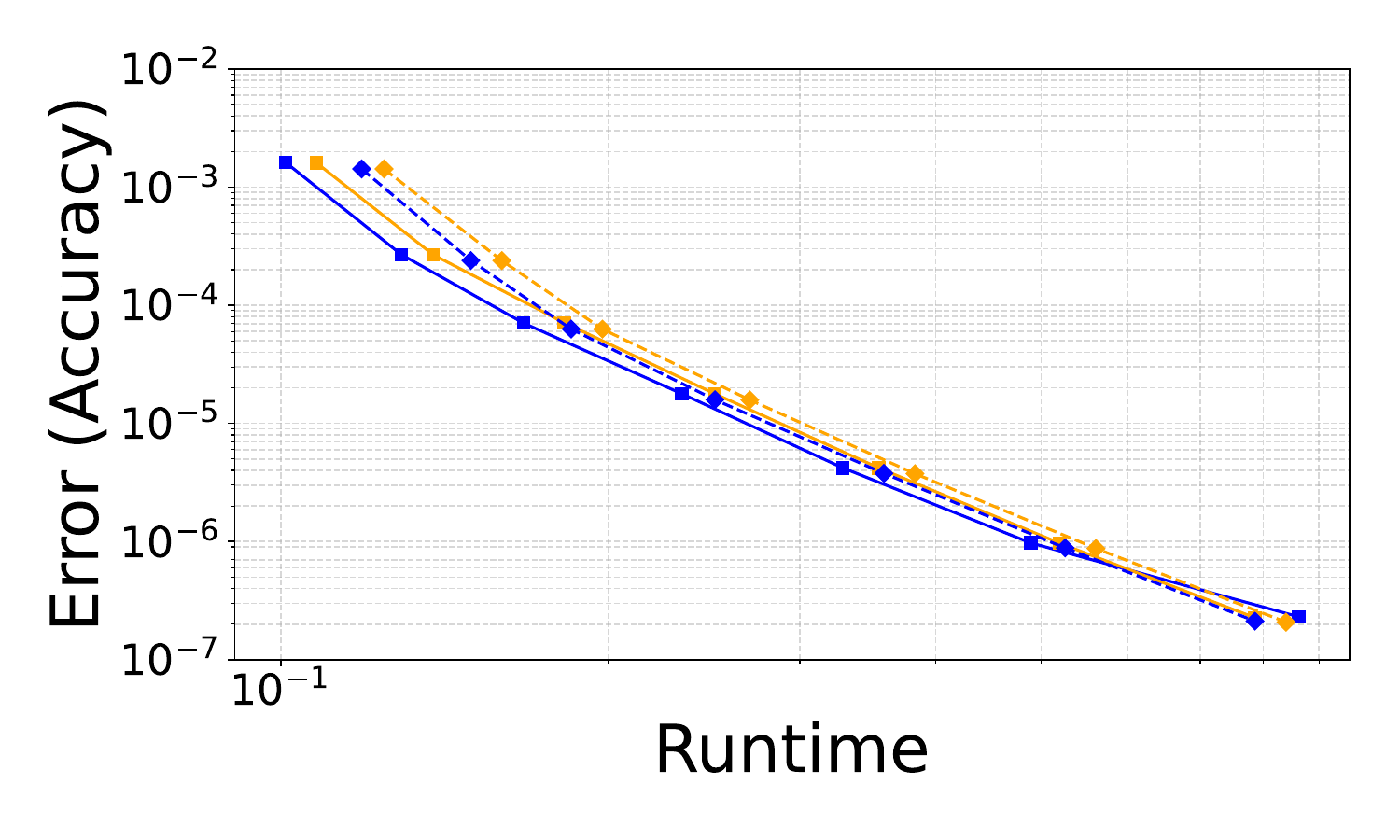}
        \caption{$\nu_{v_\|}=1.0$}
        \label{fig:error_runtime2}
    \end{subfigure}
    \hfill
    \begin{subfigure}[b]{0.305\textwidth}
        \includegraphics[trim={120 0 0 0}, clip, width=\textwidth]{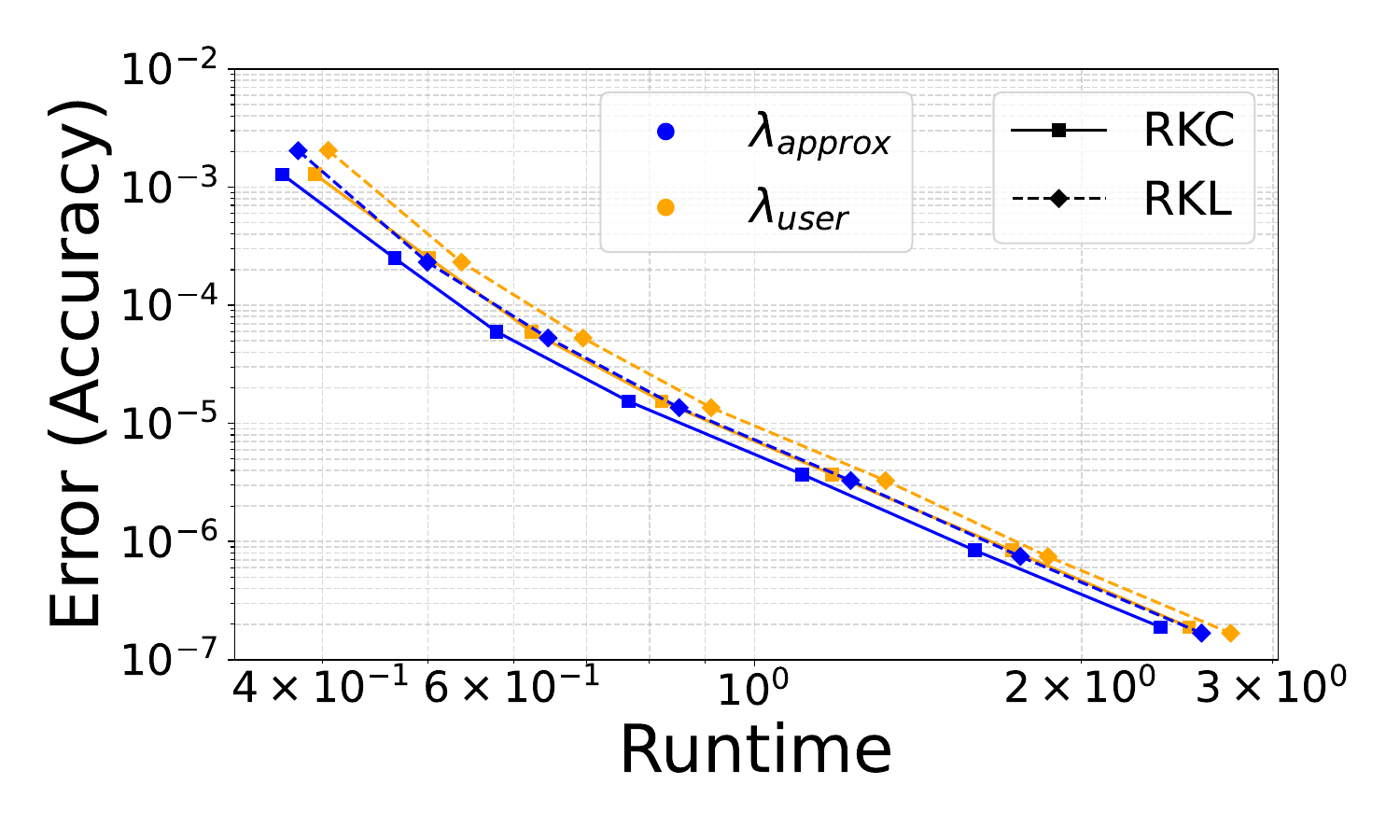}
        \caption{$\nu_{v_\|}=10.0$}
        \label{fig:error_runtime3}
    \end{subfigure}

    \caption{Computational efficiency for varying dominant eigenvalue estimation types as $\nu_{v_\|}$ increases}
    \label{fig:error_runtime}
\end{figure}

Finally, we compare results obtained with $\lambda_{user}$ and $\lambda_{approx}$. In separate computations, we noticed that the user provided dominant eigenvalue estimation, $\lambda_{user}=\overline{D}_{v_{\|}}[4/(\Delta v_{\|})]^2$ where $\overline{D}_{v_{\|}}$ stands for the cell average of the diffusion coefficient $D_{v_{\|}}$ (Eqn. \eqref{eq:diffusion_coefficients}), is approximately $16\%$ more conservative than the values estimated by SUNDIALS. Since STS methods use this estimate to decide the number of stages to take, a conservative dominant eigenvalue estimate results in a larger stability region, at the expense of additional computational effort.  We thus anticipate that a sharper estimate of the dominant eigenvalue will result in improved computational efficiency, so long as this estimate is sufficiently large as to guarantee linear stability.  Figure \ref{fig:error_runtime} confirms this hypothesis, showing that for both RKL and RKC, the simulations that use the estimated $\lambda_{approx}$ from Section \ref{sec:dominant-eigenvalue} are slightly more efficient than those that use the analytical estimate $\lambda_{user}$.

\section{Conclusions}
\label{sec:conclusions}

In this work, we have investigated the performance of diagonally implicit Runge–Kutta (DIRK) methods, explicit strong-stability-preserving (SSP) Runge–Kutta methods, and super-time-stepping (STS) methods for solving diffusion equations of relevance to gyrokinetic plasma physics. Finite-difference discretizations were studied using DIRK, SSP, and STS schemes, while discontinuous Galerkin (DG) discretizations were examined using SSP and STS methods. These approaches were implemented through coupled use of the SUNDIALS and Gkeyll frameworks, leveraging SUNDIALS' support for multiple time integration algorithms, adaptive time stepping, and dominant eigenvalue estimation, and leveraging Gkeyll's discontinuous Galerkin discretization infrastructure. Our numerical experiments demonstrate that STS methods effectively handle the stiffness associated with diffusion terms while maintaining computational efficiency, outperforming both Jacobian-free DIRK methods and explicit SSP methods for the regimes of interest.  In particular, STS methods — especially when combined with SUNDIALS-based dominant eigenvalue estimates — exhibit superior runtime performance and stability compared to SSP methods, with advantages becoming more pronounced at larger diffusion coefficients.

We further examined different error-norm formulations for adaptive time stepping and found that a cell-wise norm improves efficiency and significantly reduce step rejections, with zero rejections observed for STS methods, while the classical component-wise norm can lead to step rejection rates of up to $26\%$. Overall, these results indicate that STS methods provide a robust and efficient approach for explicit time integration in stiff diffusion-dominated problems, offering a practical alternative to fully implicit methods without incurring their associated algorithmic and implementation complexity.

Much work remains to bring these methods to bear on large-scale gyrokinetics simulations of magnetic confinement fusion devices.  Next steps include incorporation of these STS methods within multi-physics time-stepping algorithms, including operator-splitting and methods like PIROCK \cite{abdullePIROCKSwissknifePartitioned2013a}, that naturally combine STS and Runge--Kutta methods to separately tackle advection and diffusion operators present in our target application \eqref{eq:gkeyll_model}.

\section*{Acknowledgments}

We thank Ammar Hakim (PPPL) for his insight into solution of PDEs with DG methods and for sharing his experience with STS methods. This work was supported by the U.S. Department of Energy (DOE), Office of Science, Office of Advanced Scientific Computing Research, Scientific Discovery through Advanced Computing (SciDAC) Program through the Computational Evaluation and Design of Actuators for Core-Edge Integration (CEDA) project, and through the Frameworks, Algorithms and Software Technologies for Mathematics (FASTMath) Institute. This work was also supported by the DOE CEDA and Distinguished Scientist programs via contract DE-AC02-09CH11466 for the Princeton Plasma Physics Laboratory.

\bibliographystyle{unsrt}
\bibliography{references}

\end{document}